\documentclass[graybox]{svmult}

\usepackage{type1cm}        
\usepackage{nicefrac}
\usepackage{makeidx}         
\usepackage{graphicx}        
\usepackage{multicol}        
\usepackage[bottom]{footmisc}

\usepackage{newtxtext}       %
\usepackage[varvw]{newtxmath}       

\usepackage{hyperref}
\usepackage{cprotect}

\usepackage{footmisc}

\usepackage{amscd}

\usepackage{enumerate}

\newcommand{\C}{{\mathbb C}}

\newcommand{\N}{{\mathbb N}}

\newcommand{\cA}{{\mathcal{A}}}
\newcommand{\cB}{{\mathcal{B}}}
\newcommand{\B}{{\mathcal{B}}}
\newcommand{\maB}{{\mathcal{B}}}
\newcommand{\cC}{{\mathcal{C}}}

\newcommand{\cH}{{\mathcal{H}}}

\newcommand{\cP}{{\mathcal{P}}}
\newcommand{\cS}{{\mathcal{S}}}
\newcommand{\cZ}{{\mathcal{Z}}}

\newcommand{\eps}{{\varepsilon}}

\newcommand{\la}{\langle}
\newcommand{\ran}{\rangle}

\def\eps{\varepsilon}

\newcommand{\supp}{\operatorname{supp}}

\newcommand{\Aut}{\operatorname{Aut}}

\newcommand{\ds}{\displaystyle}

\newcommand{\alg}{\operatorname{alg}}

\newcommand{\be}{\begin{enumerate}}
\newcommand{\ee}{\end{enumerate}}

\begin{document}

\title*{$M$-Local type conditions for the $C^*$-crossed product
 and local trajectories}
  \titlerunning{ $M$-Local type conditions  and and local trajectories}

\author{M. Am\'elia Bastos, \\ Catarina C. Carvalho and\\ Manuel G. Dias}
\institute{M. Am\'elia Bastos \at  Instituto Superior T\'{e}cnico, University of Lisbon, Av. Rovisco
  Pais, 1049-001 Lisbon, Portugal, \email{amelia.bastos@tecnico.ulisboa.pt}
\and Catarina C. Carvalho \at  Centro de Análise Matemática, Geometria e Sistemas Din\^amicos, Instituto Superior T\'{e}cnico, University of Lisbon, Av. Rovisco
  Pais, 1049-001 Lisbon, Portugal, \email{catarinaccarvalho@tecnico.ulisboa.pt}
  \and Manuel G. Dias \at Vrije Universiteit Brussels, Dept. Mathematics, Pleinlaan 2, 1050 Brussels, Belgium,
    \email{manuel.g.dias@tecnico.ulisboa.pt}}

\maketitle



\abstract*{
The local trajectories method  establishes invertibility in  algebras  $\mathcal{B}= \alg(\mathcal{A}, U_G)$, for a unital $C^*$-algebra $\mathcal{A}$ with a non-trivial center, and a unitary group $U_g$, $g\in G$, with $G$ a discrete group, assuming that $G$ is amenable and the action  $a\mapsto U_gaU_g^*$ is topologically free. It is applicable in particular to $C^*$-algebras associated with convolution type operators with amenable groups of shifts. 
 We introduce here an $M$-local type condition that allows to establish an isomorphism between $\cB$ and a $C^*$-crossed product, which is fundamental for the local trajectories method to work. 
 We replace amenability of $G$ by the more general condition that action is amenable. The influence of the structure of the fixed points of the group action is analysed and a condition is introduced that applies when the action is not topologically free. 
If $\mathcal{A}$ is commutative, the referred conditions are related to the subalgebra $\alg(U_G)$ yielding, in particular, a sufficient condition that depends essentially on $U_G$.
It is shown that in $\pi(\mathcal{B})= \alg(\pi(\mathcal{A}), \pi(U_G))$, with $\pi$ the local trajectories representation,  the $M$-local type condition is verified, which allows establishing the isomorphism essential for the local trajectories method. 
}

\abstract{
The local trajectories method  establishes invertibility in  algebras  $\mathcal{B}= \alg(\mathcal{A}, U_G)$, for a unital $C^*$-algebra $\mathcal{A}$ with a non-trivial center, and a unitary group $U_g$, $g\in G$, with $G$ a discrete group, assuming that $G$ is amenable and the action  $a\mapsto U_gaU_g^*$ is topologically free. It is applicable in particular to $C^*$-algebras associated with convolution type operators with amenable groups of shifts. 
 We introduce here an $M$-local type condition that allows to establish an isomorphism between $\cB$ and a $C^*$-crossed product, which is fundamental for the local trajectories method to work. 
 We replace amenability of $G$ by the more general condition that action is amenable. The influence of the structure of the fixed points of the group action is analysed and a condition is introduced that applies when the action is not topologically free. 
If $\mathcal{A}$ is commutative, the referred conditions are related to the subalgebra $\alg(U_G)$ yielding, in particular, a sufficient condition that depends essentially on $U_G$.
It is shown that in $\pi(\mathcal{B})= \alg(\pi(\mathcal{A}), \pi(U_G))$, with $\pi$ the local trajectories representation,  the $M$-local type condition is verified, which allows establishing the isomorphism essential for the local trajectories method. 
}


\section*{Introduction}

The study of invertibility criteria in algebras of operators plays an important role in operator theory, with wide applications in many areas.
 One approach that has been fruitful is to use suitable families of representations, so that to reduce to invertibility of the representatives in 'nicer', so-called local, algebras. The  Allan-Douglas local principle applies to $C^*$-algebras $\cA$ with a non-trivial center, in that if $\cZ$ is a central subalgebra then using the isomorphism $\cZ \cong C(M)$ given by Gelfand theory, with $M$ the maximal ideal space of $\cZ$, we have local algebras 
 $$\cA_m = \cA / J_m, \quad \pi_m : \cA\to \cA_m, \quad m\in M,$$
 where $J_m$ is the closed *-ideal of $\cA$ generated by $m$, $\pi_m$ is the quotient map, and in this case, $a\in\cA$ is invertible  in $\cA$ if, and only if, $\pi_m(a)$ is invertible in $\cA_m$, for all $m\in M$. 
 
  In this paper, we are interested in the case when we have a unital $C^*$-algebra of operators $\cA$, a discrete group $G$ and a group of unitary operators $U_G$ defining an action of $G$ on $\cA$, that is, if we have a $C^*$-dynamical system. One is then bound to study invertibility in the algebra of operators generated by $\cA$ and $U_g$, $g\in G$, 
  $$\cB := \alg (\cA, U_G), \quad \mbox{with  }  \cA\subset B(H), \quad U: G\to B(H) \mbox{  unitary,}$$
with $B(H)$ be the $C^*$-algebra of all bounded linear operators acting on some Hilbert space $H$. The issue is that even when $\cA$ is commutative, $\cB$ typically  has a trivial center, so local principles do not apply directly.  

The study 
 of  $C^*$-algebras of operators associated with $C^*$-dynamical systems was
  developed by Antonevich when the initial $C^*$- algebra is commutative and the group is subexponential or admissible (see \cite{Anton_paper, Anton_book}). 
The extension of the Allan-Douglas local principle  to  $C^*$-algebras 
 induced by an action of a discrete amenable group 
   was developed by Y. Karlovich  \cite{karlovich88, Karlovich_OTAA17}, relying on a given arbitrary central subalgebra of $\cA$ and representations of $C^*$-crossed product algebras, 
   originating what we call the \emph{local trajectories method} (see below for details). 
      An alternative approach was also developed by Antonevich, Lebedev, Brenner (see \cite{AntonLebedev} and references therein). 
   
   The local trajectories method gives a powerful machinery for studying invertibility and Frehdolmness in $C^*$-algebras of nonlocal type operators with discontinuous data. A first example of this application can be found in the paper by Yuri Karlovich \cite{karlovich88-2} in which he analysed the Fredholm theory of convolution type operators with discrete groups of displacements and coefficients admitting discontinuities of semi almost periodic type. 
It has also been applied in the Fredholm analysis of $C^*$-algebras with amenable groups of shifts, algebras of convolution type operators with oscillating coefficients, algebras of singular integral operators with piecewise quasi-continuous and semi-almost periodic coefficients.    Examples of these applications can be found  for instance in 
   \cite{BFK18, BFK_PQC_2019, BFK_Spectral07, BKarS_8, karlovich88-2, ShvyLatu}.

More precisely,
 the local trajectories method is applicable to study invertibility in $C^*$-algebras of the form $\cB=\alg (\cA,U_G )$, where $\cA\subset B(H)$ is a $C^*$-algebra with non-trivial center and $U_G:= U(G)$, where $U: G\to B(H)$  is a unitary representation of an amenable discrete group $G$. We assume that we have an action $\alpha: G\to \Aut(\cA)$ given by $\alpha_g(a)= U_g a U_g^*$, so we have a dynamical system $(\cA, G, \alpha)$, with $(id, U)$ a covariant representation on $H$, and we can consider the  crossed product algebra $\cA\rtimes_\alpha G$. We assume that $G$ also acts on a given central subalgebra $\cZ\cong C(M)\subset \cA$, so that we have an induced action of $G$ on the space $M$. 
 
 The idea behind  the local trajectories method is to consider first localization in $\cA$  as in Allan-Douglas and then consider the associated regular representations of $\cA\rtimes_\alpha G$ on $\ell^2$-spaces. If we let $\Omega$ be the set of orbits of the action on $M$,  the for points in the same orbit the respective regular representations are unitarily equivalent, so the local trajectories representations reduce to a family of  regular representations $\{\pi_\omega\}_{\omega\in \Omega}$ of the crossed product algebra $\cA\rtimes_\alpha G$. (We review these definitions in Section  \ref{s.localtraj}.)
  
   In order to study invertibility in $\cB$ through this family of representations,  it is therefore a fundamental  step to establish an isomorphism between the algebra $\cB$ and the crossed product $\cA\rtimes_\alpha G$. Conditions for this isomorphism, as well as for the faithfulness of the local trajectories family, typically assume that the group $G$ is amenable. The approaches in \cite{AntonLebedev, Anton_book, karlovich88, Karlovich_OTAA17} rely on proving  suitable isomorphism theorems giving sufficient conditions namely through the set of fixed points and the crucial notion of a \emph{topologically free} action (see also \cite{ArchboldSpielberg}).

In general, an action is said to be topologically free if considering the dual action  induced on the spectrum $\hat{\cA}$ of $\cA$, the sets of fixed points by finite, non-trivial, subsets of $G$ have empty interior (one can also consider the dual action on the space of primitive ideals as in \cite{AntonLebedev}). In the setting of local trajectories, following \cite{Karlovich_OTAA17},  the notion of being topologically free relies on the topology of the pure state space.
If we let $\cP_\cA$  and $\cP_\cZ$ be  the classes of pure states of $\cA$ and $\cZ$, then since $\cZ\cong C(M)$, we have $\cP_\cZ \cong \hat \cZ \cong M$, and 
 since $\cZ$ is central, there is a well-defined, surjective,  restriction map
$$\psi: \cP_\cA \to \cP_\cZ \cong \hat \cZ \cong M.$$
Then   the action is said to be \emph{topologically free} if for any finite set $G_0\subset G$,  and for any non-empty open set $W\subset \cP_\cA$,  there exists $m\in \psi(W)$ with $\beta_g(m)\neq m$, for all $g\neq e \in G_0$, with  $\beta : G\times M \to M$ the induced action on $M$; this is referred to as condition (A3) (see Section  \ref{s.localtraj}). 

 The main result in \cite{Karlovich_OTAA17} (see also Theorem  \ref{thm.mainlocaltrajec}) can then be written as:

\emph{Local Trajectories Method: If $G$ is amenable and the action is topologically free, then $\cB\cong \cA\rtimes_\alpha G$ and the family $\{\pi_\omega\}_{\omega\in \Omega}$ is faithful, so that $b\in \cB$ is invertible if, and only if, $\oplus_{\omega} \pi_\omega(b)$  is invertible if and only if, $\pi_\omega(b)$ is invertible, for all orbits $\omega$, and $\sup_{\omega\in\Omega}\|\pi_\omega^{-1}(b)\|<\infty$.}

 If the action is not topologically free, the situation is typically much harder to analyse, but in some situations there are still methods to reduce to this case, see  \cite{BFK_Spectral07, Karlovich_OTAA17} (and references therein), and also  \cite{AntonLebedev, Anton_book}.

The purpose of the present article is to explore the relation between the notion of topologically freeness and the isomorphism with the crossed product algebra, as well as alternative conditions for the local trajectories method to hold in order to better understand its domain of applicability. 

We do this by going back to a basic global condition, known to be equivalent to the existence of an isomorphism $\cB\cong   A\rtimes_\alpha G$ - referred to as \emph{condition (B0)} - 
 and establishing similar conditions but of \emph{$M$-local type}, that is, on open sets of the space $M$ of maximal ideals in $\cZ$ - referred to as \emph{condition (B1)}. This basically amounts to considering 'local' elements of the form $zb$, where $b\in \cB$ and $z\in \cZ=C(M)$ is a bump function with support in $V$ (the precise definitions are in Section \ref{ss.gencond}).

We note that a necessary condition for condition (B0) to yield the isomorphism of $\cB$ with the crossed product is that the full and reduced crossed product algebras coincide. In previous works, this identification typically comes from assuming that the group $G$ is amenable. Here this condition is replaced by the more general notion of \emph{amenability of the action}  $\alpha: G \to \Aut(\cA)$ (see for instance \cite{ADelarocheAmean, BrownOzawa_book}), which we assume from Section \ref{s.cond_isom}. That is, instead of assuming that the group is 'nice', we assume that it acts 'nicely', and we still have $\cA\rtimes_\alpha G \cong \cA\rtimes^r_\alpha G$.

 We show that  $M$-localization  applies to a  class of algebras $\cA$ where open sets of $\cP_\cA$ always determine an open set in $M$, in that the restriction map $\psi$ defined above is open. We say that such an algebra $\cA$ is \emph{$M$-localizable}, or satisfies condition (C). 
 This  class includes all commutative algebras and algebras of matrices of continuous functions, as well as the class of algebras  considered in \cite{Anton_book}.

 The point in using $M$-local type conditions is that we can split between open sets $V\subset M$ where the action is 'well-behaved', in that there are points in $V$ that are not fixed  by a finite $G_0\subset G$, and the open sets $V$ that are fixed by some  non-trivial finite $G_0\subset G$ such that $(\beta_g)_{|V}=id_{|V}$, $g\in G_0$
 (see Section  \ref{ss.fixed_points}). This is made possible by Lemma \ref{lem_fixedpts} and Proposition \ref{prop.B2prev},  which imply on one hand that if the action is topologically free then the $M$-local condition (B1) is automatically satisfied, so that we have  an isomorphism with the crossed product, and on the other hand, that in the general case, to check (B1) it suffices to consider the class of open sets where the action of some finite subset is trivial - referred to as \emph{condition (B2)}. 
 As a result, in Corollary \ref{cor.B2equivB1}, we obtain an isomorphism between $\cB$ and $\cA\rtimes_\alpha G$ as long as we replace the assumption of the action being topologically free
  by our $M$-local condition being satisfied on such sets, should they exist, condition (B2).
 
 We use condition (B2) to show, in Examples \ref{example1} and \ref{example2},  that in case the action is not topologically free, then the fact that $\cB$ is isomorphic to the crossed product may depend on the way $\cA$ and $\alg(U_G)$ ’sit’ inside $B(H)$ and how they interact, which cannot happen in the presence of topological freeness, due to the isomorphim theorems in \cite{AntonLebedev, ArchboldSpielberg, Karlovich_OTAA17}.

 This work then explores the local type conditions in two directions. The first one, developed in Section \ref{ss.subalgebras_comm}, applies directly to commutative algebras, and has to do with a $C^*$-subalgebra $\cA'$ that is $M$-locally 'arbitrarily close' to $\cA$ and where $G$ also acts; in this case, assuming the generated algebras are also $M$-locally 'arbitrarily close', then the local condition (B2) needs only be verified on $\cA'$ and $\cB'$. We say that $\cA'$ and $\cB'$ are \emph{$M$-locally dense in $\cA$ and $\cB$.}
 
 It turns out that when the algebra $\cA$ is commutative,  we see that this is the case with $\cA'=\C$ and $ \cB'=\alg(U_G)$, 
 so we  obtain a  $M$-local condition in $ \alg(U_G)$  sufficient for our local condition (B2) to be verified in $\cB$, and for the isomorphism with the crossed product algebra; this is Theorem \ref{prop.comm_red}.

 Moreover, under an additional non-degeneracy condition, we can use (B2) to get rid of the $M$-local element, and we find a global condition on $ \alg(U_G)$, namely that $ \alg(U_G)\cong C^*(G)$ the group algebra, which guarantees (B2) and  the isomorphism with the crossed product algebra. We point out that this condition depends only on $U_G$, not on localization, and effectively substitutes the requirement of topological freeness in this situation, see Theorems \ref{thm.condscomm2} and \ref{thm.condslocadense2.localtrajec}.

Another application of the $M$-local conditions (B1) and (B2) is to guarantee conditions for applicability of the method of local trajectories for $M$-localizable algebras, that is,  satisfying (C). As it is pointed out in Section \ref{s.localtraj}, in general, the local trajectories method is applicable to $\cB$, and yields a faithful family, if there are isomorphisms
$$ \cB \cong A\rtimes_\alpha G, \qquad \pi(\cB)\cong \pi(\cA) \rtimes_{\alpha'} G$$
where $\pi=\oplus \pi_\omega$ and $\{\pi_\omega\}_{\omega \in \Omega}$ is the local trajectories family, $\pi_\omega$ regular representations.
For the second isomorphism,   we have representations on $C^*$-algebras  of the form $B(\ell^2(G,H_\omega))$  and we can explore their norm properties to show that 
 the $M$-local type condition (B1) is in fact always verified 
   in $\pi(\cB)= \alg(\pi(A), U'_G)$. 
   We conclude that if $\cA$ is $M$-localizable, then the conditions obtained before to ensure the isomorphism $ \cB \cong A\rtimes_\alpha G$ are also sufficient  for the method of local trajectories to work on $\cB$; this leads to Theorems \ref{thm.extracond}, \ref{thm.condscomm2.localtraj} and \ref{thm.condslocadense2.localtrajec}.

 In this article, we chose to formulate our results for a 'concrete' algebra of operators $\cA\subset B(H)$, for some fixed Hilbert space $H$,  together with an action of a discrete group on $\cA$, and a group of unitary operators $U_G\subset B(H)$. The condition that the action satisfies $\alpha_g(a)= U_gaU_g^*$ is nothing more than the pair $(id, U)$ being a covariant representation for $(\cA, G, \alpha)$. 
 We point  out that in fact our results also apply  to the  setting when we start with an 'abstract' unital $C^*$-algebra $\cA$ and consider an arbitrary faithful representation $\phi: \cA\to B(H_\phi)$ on some Hilbert space $H_\phi$, and a covariant representation $(\phi, U)$, with  $U_G\subset B(H_\phi)$  unitary operators. Defining $\cB_{\phi, U}= \alg (\cA, U_G)$, conditions (B0), (B1) and (B2) can be easily written  (see Remark \ref{rmk.gencovrepB0}) in a way such that they yield an isomorphism   $\B_{\phi,U}\cong  \cA\rtimes_\alpha G $, and similar conditions can be given such that the locally trajectories method works on $\cB_{\phi, U}$. 


We now give an outline of the paper. In Section \ref{s.prelim}, we  review some objects and concepts needed throughout the paper.
 We also establish the setting we will work on, 
 and define the relevant $C^*$-algebras associated with our structures. 

  In Section \ref{s.localtraj}, we present the local trajectories method.
    We follow the lines of  \cite{Karlovich_OTAA17}, while focusing on pinpointing the main steps and in particular, in the role of the isomorphisms with crossed product algebras. 

In Section \ref{s.cond_isom}, we explore conditions to guarantee the isomorphism with the crossed product,  using what we call $M$-localization, that is, reducing to open subsets of $M$, and assuming that the action is amenable. 
 
  We start, in Subsection  \ref{ss.gencond}, by introducing the classical condition to guarantee the isomorphism of $\cB$ with the respective crossed product algebra, called here condition (B0). We then give a similar condition of $M$-local type,  on open sets of the space of the maximal ideals of $\cZ$, called condition (B1). Under an assumption on $\cA$, condition (C) that $\cA$ is $M$-localizable, we prove the equivalence between these two notions, Proposition \ref{localB1}.

In Subsection \ref{ss.fixed_points} we see how this localized condition can be used together with the structure of the fixed point set, and the notion of topologically free action, to guarantee the isomophism with the crossed product algebra, and arrive at condition (B2) that 
  only involves open sets of $M$ fixed by a finite subset of $G$ (Proposition \ref{prop.B2prev}).

 In Subection \ref{ss.subalgebras_comm}, we  show that in the commutative case, we need only check conditions for isomorphism on $\alg(U_G)$ and the group algebra $C^*(G)$, depending only on $G$.  We prove first that  this algebra is what we call $M$-locally dense in $\cB$, and abstract these results to the general case.
 
  Then in Section \ref{ss.localtraj_cond}, we tackle the conditions for applicablity of the local trajectories method, applying our results
   to the local trajectories representation $\pi= \oplus \pi_\omega$. 
  We show that the image $\pi(\cB)$ of $\cB$ 
   is always isomorphic to a crossed product, proving that if $\cA$ is $M$-localizable, the local trajectories method works on $\cB$ as long as there is an isomorphism of $\cB$ with the crossed product. In the commutative case, we give a sufficient condition depending only on $U_G$.

\textbf{Acknowledgments}: The authors would like to thank an anonymous referee for carefully reading the manuscript and contributing with various comments and suggestions that helped to improve the paper.

\section{Preliminaries}\label{s.prelim}

Throughout the paper, $\cA$  will always denote a unital $C^*$-algebra. 
 We  review in this section some concepts and  results on $C^*$-algebras and representation theory that will be needed in the paper. 
 For general references, see for instance \cite{Dixmier, Kadison-Ring_book, Pedersen_book}.  

\subsection{Representations and states}\label{ss.rep_states}

By a \emph{representation} of a $C^*$-algebra $\cA$
on a Hilbert space $\cH_\pi$ we always mean a non-degenerate $*$-homomorphism $$\pi : \cA \to
B(H_\pi),$$
with $B(H_\pi)$ the $C^*$-algebra of bounded linear operators on $H_\pi$. 
 Non-degenerate in the unital case means that $\pi(1_\cA)=I$. 
We say that $\pi$ is \emph{irreducible} if its only invariant subspaces are trivial. An injective representation is called \emph{faithful}.

  A \emph{state} in $\cA$ is a positive linear functional $\mu$ on $\cA$ with  $\|\mu\| =1$, or, equivalently, $\mu(1_\cA)=1$. The state space $\mathcal{S}_\cA\subset \cA^*$ is convex and compact in the weak-$*$ topology, and as such, 
  has extreme points that are called \emph{pure states}. We denote the pure state space of $\cA$ by $\cP_\cA$, and always endow $\mathcal{S}_\cA$ and $\cP_\cA$ with the weak-$*$ topology, that is the topology of pointwise convergence.   
  
  States are related to representations in a fundamental way: $\pi: \cA\to B(H)$ is an irreducible representation on a Hilbert space $H$, and $x$ is a unit vector in $H$, then
 $$\mu(a):= \la \pi(a) x, x\ran$$
 is a pure state, and, by the GNS construction, any pure state defines a irreducible representation $\pi_\mu$ satisfying the above.

By the Gelfand-Naimark theorem, for every $C^*$-algebra $\cA$ there exists a faithful representation $\pi: \cA\to B(H)$
 given by the direct sum 
$$\pi= \oplus \pi_{\nu}, \quad \pi_{\nu}: \cA \to B(H_{\nu})$$
where $\pi_{\nu}$ are{ irreducible} representations associated to pure states. 

 We will make extensive use of the following properties (see \cite{Kadison-Ring_book}, Proposition 4.3.1 and Theorems 4.3.8 and  4.3.14):
 
 \begin{enumerate}[(i)]
 \item Pure states separate points in $\cA$: if $\mu(a)=0$, for all $\mu\in \cP_\cA$ then $a=0$. 
 
 \item  For any state $\mu$ on $\cA$, and $a, b\in \cA$, we have the Cauchy-Schwarz inequality for states:
 $$|\mu(a^*b)|^2\leq \mu(a^*a)\mu(b^*b).$$
{In particular, since $\mu(1_\cA)=1$, we get $|\mu(a)|\leq \sqrt{\mu(a^*a)}\leq \|a\|$.}
 
\item  {If $a\in \cA$ is normal, that is, if $a^*a=aa^*$, then there exists a pure state $\mu\in \cP_\cA$ such that
 $$\|a\|= |\mu(a)|.$$
 } 
 In general, for $a\in \cA$, 
 $$\|a\|= \max_{\mu\in \cP_\cA} \sqrt{{\mu(a^*a)}}.$$  
 In particular, if $a$ is a positive element, then $\|a\|= \max_{\mu\in \cP_\cA} \mu(a)$. 
 
 \item If $\cZ$ is a central subalgebra of $\cA$, then 
 for any $\mu\in \cP_\cA$, we have 
 $$\mu(za)= \mu(z)\mu(a), \quad \mbox{ for } a\in \cA, z\in \cZ.$$
 It follows that $\mu_{|\cZ}$ is a multiplicative linear functional. 
  \end{enumerate}
 It follows from (iii) and the GNS construction that for any $a\in \cA$ there exists an irreducible representation $\phi$ of $\cA$ such that 
 \begin{equation*}\label{eq.irredrep_norming}
 \|a\|= \|\phi(a)\|. 
\end{equation*}
We shall also need results on extension of states. Let $\cZ$ be a closed $C^*$-subalgebra of $\cA$, containing the identity. Then any state on $\cZ$ can be extended to a state on $\cA$
  (see \cite{Kadison-Ring_book}, Theorem 4.3.13); for each state in $\mathcal{S}_\cZ$, the set of its extensions to $\mathcal{S}_\cA$ is weak-$*$compact and convex. Moreover, pure states can be extended to pure states. 
  {On the other hand, if $\cZ$ is a \emph{central} subalgebra of $\cA$, pure states restrict to pure states (by (iv) and Proposition 4.4.1 in \cite{Kadison-Ring_book})}, so restriction yields then a surjective map 
  \begin{equation}\label{eq.def.PAtoM}
 \psi: \cP_\cA \to \cP_\cZ, \qquad\mu \mapsto \mu_{|\cZ},
 \end{equation}
 which is moreover continuous, since $\cP_\cA$, $\cP_\cZ$ have the weak-$*$ topology.
Since $\cZ$ is a central subalgebra, then from (iv) above, 
 $\mu_{|\cZ}$ is a multiplicative linear functional and $\cP_\cZ = \hat \cZ$, the character space of $\cZ$. Moreover, by Gelfand's theorem, 
  since $\cZ$ is commutative, $\hat \cZ\cong M$, the space of maximal ideals of $\cZ$, and 
 the Gelfand transform yields an isomorphism $\cZ\cong C(M)$.

%

\subsection{Algebras associated to an unitary action}

In what follows, we let $\cA$ be a unital  $C^*$-algebra and $G$ be a discrete group. For details on the constructions below, see for instance \cite{BrownOzawa_book, Pedersen_book, dwilliams_book} and references therein.

An \emph{action} of $G$ on $\cA$ is a homomorphism $\alpha : G \to  \Aut(\cA)$, where $ \Aut(\cA)$ is the group of $*$-automorphisms of $\cA$. Given such an action we call  $(\cA, G, \alpha)$ a \emph{$C^*$-dynamical system}. 

Given a $C^*$-dynamical system $(\cA, G, \alpha)$, we denote by  $C_c(G,\cA)$  the linear space of finitely 
 supported functions in $G$,
$$C_{c}(G,\cA) =\{f: G\to \cA\, | \,  f(s)=0 , s \notin G_{0} \mbox{  finite }\}.$$
We use the action $\alpha$ to define a $\alpha$-twisted convolution product on $C_{c}(G,\cA) $, as well as an involution:
$$(f*g)(s):= \sum_{t\in G} f(t)\alpha_{t}(g(t^{-1}s)), \qquad f^{*}(s):= \alpha_{s}(f^{*}(s^{-1})).$$
We call  the $*$-algebra $C_{c}(G,\cA)$ the 
\emph {convolution algebra} of $(\cA, G, \alpha)$.  

If $\cA$ is commutative, $\cA\cong C(M)$, with $M$ some compact Hausdorff space, then  $C^*$- dynamical systems are in one-to-one correspondence with group actions on $M$: if $G$ acts on compact space $M$, $G\times M\to M, (g,m) \to g\cdot m$,  
  then
$(C(M), G, \alpha)$ is a $C^*$- dynamical system with
$$\alpha_{s}(f)(m)=f(s^{-1} \cdot m), \quad m\in M, s\in G.$$
Conversely, if $\cA= C(M)$ is a commutative algebra, then $G$ acts on $M$ and $(G,M)$ is a transformation group.

  We are interested in studying invertibility in $C^*$-algebras  associated to dynamical systems. 
  
 \begin{definition}\label{def.algB}
 Assume that $\cA\subset B(H)$, for some Hilbert space $H$, and let $U: G\to B(H), g\mapsto U_g$ be a unitary representation. We denote by $$\B:= \alg(\cA, U_{G})$$ the $C^*$-subalgebra of $B(H)$ generated by  $\cA$ and $U_G=\{U_g, g\in G\}$. Assume also that $U_{g}a U^{*}_{g}$ is a $*$-automorphism of $\cA$, for all $g\in G$, that is, that we have an action
 $$\alpha: G \to \Aut(\cA), \qquad \alpha_{g}(a):= U_{g}a U^{*}_{g},$$
  then   $\B$ is the closure in $B(H)$ of the $*$-subalgebra
 \begin{equation*}
 \B_0:= \left\{ \sum_{g\in G_0} a_g U_g: a_g\in \cA, G_0\subset G \mbox{ finite }\right\}=  \left\{ \sum_{g\in G} a(g) U_g: a\in C_c(G,\cA)\right\}.
 \end{equation*}

 \end{definition}

Given an arbitrary dynamical system  $(\cA, G, \alpha)$, there is always a universal object that encodes both the original $C^*$-algebra $\cA$ and the group action. 

\begin{definition}
The \emph{crossed product algebra} $\cA  \rtimes_{\alpha} G$ is the completion of $C_{c}(G,\cA) $ with respect to the universal norm 
 $$\|f\|_{u} = \sup_{\pi} \|\pi(f) \|,$$
where $\pi$ ranges over all  $*$-homomorphisms $\pi: C_c(G, \cA) \rightarrow$ ${B}({H})$, with $H$ a Hilbert space. 
\end{definition}

When $\cA=\C$, the crossed product algebra yields the group algebra $C^*(G)$.

If $\cB$ is as in Definition \ref{def.algB}, one of our goals is to discuss conditions under which $\cB \cong \cA  \rtimes_{\alpha} G$ (see Section \ref{s.cond_isom}). One can check that there is always a $*$-homomorphism $\Phi: C_{c}(G,\cA)\to \B_{0}$, surjective, given by
\begin{equation}\label{eq.Phi}
\Phi(f) = \sum_{g\in G} f(g) U_{g},
\end{equation}
that  
  extends to a surjection $\Phi: \cA  \rtimes_{\alpha}G \to \B$, so the algebra $\cB$ is always a quotient of the crossed product algebra. 
To study representations of such algebras $\cB$, we consider first representations of the crossed product. 

Let $(\cA, G, \alpha)$ be a dynamical system. Given  a Hilbert space $H$, consider a representation $\pi: \cA\to B(H)$ and  a unitary representation $U: G\to B(H)$, then the  pair $(\pi, U)$ is said to be a  \emph{covariant representation} of $(\cA, G, \alpha)$ if 
$$ \pi(\alpha_{s}(a))= U_{s} \pi(a) U_{s}^{*}.$$
If $(\pi, U)$ is a covariant representation on some Hilbert space $H$, we can  define the so-called  \emph{integrated representation } of  $\cA \rtimes_{\alpha} G$ by
$$\pi \rtimes U : \cA \rtimes_{\alpha} G \to B(H), \qquad (\pi \rtimes U)(f) =  \sum_{t\in G}{\pi}(f(g)) U_{g},  \quad f\in C_{c}(G,\cA).$$
We have  that any covariant representation defines a $*$-representation of $C_c(G,\cA)$ and, conversely, every non-degenerate $*$-representation of $C_c(G,\cA)$ is induced by some covariant representation of $(\cA, G, \alpha)$, so that 
$$\|f\|_{u} = \sup_{(\pi,U)} \|\pi \rtimes U(f) \|,$$
where $(\pi, U)$ ranges over all covariant representations of $(\cA, G, \alpha)$. 
 Note that in Definition \ref{def.algB}, we are simply assuming that $(id, U)$ is a covariant representation, and the map  $\Phi$ in \eqref{eq.Phi} is given by $\Phi= id\rtimes U$. 

 In fact, the definition of the crossed product algebra yields the following \emph{universal property}: for every covariant representation $(\pi, U)$, there is a $*$-homomorphism $\sigma: \cA \rtimes_\alpha G \rightarrow  \alg(\pi(\cA), U_G) \subset {B}({H})$ such that
 $$
\sigma\left(f\right)=\sum_{g\in G} \pi\left(f(g)\right) U_g, \quad  \text{for all} \; f \in C_c(G, \cA).
$$
We are interested in a particular class of integrated representations that correspond to taking the \emph{left regular representations of $G$}: given  a Hilbert space $H$, one defines the representation 
$$\lambda: G\to B(\ell^{2}(G,H)), \qquad \lambda_{g}\xi(s)= \xi( g^{-1} s).$$
 Given a representation $\pi: \cA \to B(H)$  of $\cA$, define also
$$\tilde \pi : \cA \to B(\ell^{2}(G,H)), \qquad \tilde \pi (a)\xi(s)= \pi( \alpha_{s}^{-1}(a))\xi(s).$$
Then one can see that $(\tilde \pi, \lambda)$ is {covariant representation}.
\begin{definition}\label{def.redcrossed}
The \emph{regular representation  of $ \cA \rtimes_{\alpha} G$ induced by the representation $\pi$} is given by the integrated representation induced by $(\tilde \pi, \lambda)$, that is, 
\begin{equation*}
\tilde \pi \rtimes \lambda : \cA \rtimes_{\alpha} G \to  B(\ell^{2}(G,H))
\end{equation*}
such that for $f\in C_c(G,\cA)$, 
\begin{equation}\label{def.regrep}
[(\tilde \pi \rtimes \lambda)(f) \xi](g)=\sum_{s \in G} \tilde{\pi}\left(\alpha_{g^{-1}}[f(s)]\right) \xi\left(s^{-1} g\right), \quad \xi \in \ell^2(G, H).
\end{equation}
 The 
\emph{reduced crossed product algebra} 
$ \cA  \rtimes_{\alpha}^{r} G$  is the completion of $C_{c}(G,\cA) $ in the norm 
$$\|f\|_{r}= \sup_{(\tilde \pi, \lambda)} \|(\tilde\pi\rtimes \lambda)(f)\|,$$
where $(\tilde\pi\rtimes \lambda)$ ranges over all regular representations of $ \cA \rtimes_{\alpha} G$.
If $\cA=\C$ then $ \cA  \rtimes_{\alpha}^{r} G$  is the reduced group algebra $C^*_r(G)$. 
\end{definition}

In this paper, we shall work in the setting where the full and reduced crossed product algebras coincide,  so that
 $$\|f\|_u =\|f\|_{r}= \sup_{(\tilde \pi, \lambda)} \|(\tilde\pi\rtimes \lambda)(f)\|, \quad \text{with}\; (\tilde\pi\rtimes \lambda) \text{ a regular representation}.$$
 A discrete group $G$ is said to be \emph{amenable} if there exists a state $\mu$ on $\ell^\infty(G)$ which is invariant under the left translation action: for all $s\in G$ and $f\in \ell^\infty(G)$, 
 $$\mu(sf)= \mu(f).$$
 The state $\mu$ is called an invariant mean. 
  In this case, we have $C^*_r(G)\cong C^*(G)$ and  
it follows that  $\cA  \rtimes_{\alpha} G\cong \cA  \rtimes_{\alpha}^{r} G$, for any dynamical system $(A, G, \alpha)$. 
The class of amenable groups includes all compact groups, abelian groups, solvable groups 
and finitely generated groups of subexponential growth. On the other hand, if $G$ contains a copy of the free group in two generators, then $G$ is not amenable.

We will consider here a more general notion of amenability that suffices for our purposes, that of an \emph{amenable action} (see for instance \cite{BrownOzawa_book}, Section 4.3). Define a norm in $C_c(G,\cA)$ by
$$\|f\|_2:= \left\| \sum_{g\in G} f(g)f(g)^*\right\|^{1/2}, \quad f\in C_c(G,\cA).$$
Let $\mathcal{Z}(\cA)$ denote the center of $\cA$. For $s\in G$, let $\delta_s\in C_c(G,\cA)$ be such that $\delta_s(s)=1_{\cA}$ and $\delta_s(g)=0$, $g\neq s$. 
Then $(\delta_s * f)(g)= \alpha_s ( f(s^{-1}g) ),$ $g\in G$. 
\begin{definition} An action $\alpha: G \rightarrow \operatorname{Aut}(\cA)$ is \emph{amenable} if there exist finitely supported functions $x_i: G \rightarrow \cZ(\cA)\subset  \cA$, $i\in \N$, with the following properties:
\begin{enumerate}
\item  $x_i(g) \geq 0$  for all $i \in \mathbb{N}$ and $g \in G$;
\item  
$\sum_{g \in G} x_i(g)^2=1_{\cA}$  for all $i \in \mathbb{N}$;
\item $\left\|\delta_s * x_i-x_i\right\|_2 \rightarrow 0,$
 for all $s \in G$.
\end{enumerate}
\end{definition}

Given a dynamical system $(\cA, G, \alpha)$ where  the action is amenable, we always have  $\cA  \rtimes_{\alpha} G\cong \cA  \rtimes_{\alpha}^{r} G$ (Theorem 4.3.4 in \cite{BrownOzawa_book}).
 
 Moreover, every action by an amenable discrete group is amenable.
This can be seen using the following equivalent definition of amenability: a discrete group $G$ is amenable if, and only if, for any  finite set $G_0\subset G$ there exist  nonnegative unit vectors  $x_i\in \ell^2(G)$, $i\in \N$, such that $\|\lambda_sx_i-x_i\|_2\to 0$, for all $s\in G_0$ (see for instance  \cite{BrownOzawa_book}, Theorem 2.6.8). In this case, for any action $\alpha$ of an amenable discrete group $G$ on a $C^*$-algebra $\cA$,  identifying $\mathbb{C}$ with  $\mathbb{C}1_\cA\subset \cZ(\cA)$, 
we have  $(\delta_s * x_i)(g)= \alpha_s ( x_i(s^{-1}g) )= x_i(s^{-1}g)\alpha_s ( 1_\cA )= \lambda_sx_i(g)1_\cA,$ $g\in G$.
  On the other hand, there are many relevant amenable actions of non-amenable groups, such as the action of a free group on its Gromov boundary, see \cite{BrownOzawa_book}.


\section{Local trajectories method}\label{s.localtraj}

We let $\cA\subset B(H)$ be a unital $C^*$-algebra and $U: G\to B(H)$ be a unitary representation.
We review in this section the method of local trajectories, whose goal  is to establish an invertibility criterion for operators in $\mathcal{B}=\alg(\cA, U_G)$ in terms of the invertibility of local
 representatives associated to the orbits of an action of $G$ on some compact space. 
 We follow here the approach and notation in \cite{Karlovich_OTAA17} (see also \cite{Anton_book}).

From now on, we will let $\cZ$ be a central subalgebra  of $\cA$, with $1_\cA\in \cZ$. 
In this case, the Gelfand transform yields an isomorphism $\cZ\cong C(M)$, where $M $ is the compact Hausdorff space of maximal ideals of $\cZ$ or equivalently, the class of non-zero multiplicative linear functionals with the weak-$*$ topology. 
 We typically identify $z\in \cZ$ with its Gelfand transform in $C(M)$, that is, we regard $z$ as a continuous function.
 
Following the terminology in \cite{Karlovich_OTAA17}, we consider conditions (A1) and (A2):

\emph{Condition (A1): For every $g \in G$, the mapping $\alpha_g: a \mapsto U_g a U_g^*$ is a $*$-automorphism of the $C^*$-algebras $\mathcal{A}$ and $\mathcal{Z}$}.

\emph{Condition (A2): $G$ is an amenable discrete group.}

Assume condition (A1), that is, we assume that  $(id, U)$ is a covariant representation of $(\cA, G, \alpha)$ where $\alpha$ is an action on $\cA$ that maps $\cZ= C(M)$ to $\cZ=C(M)$.  
  We  have then an action $\beta$ of $G$ on $M$ given by 
$ \beta_{g}: M\to M $ %
 such that for $z \in \mathcal{Z}$,
$$
z\left(\beta_g(m)\right)=\left(\alpha_g(z)\right)(m), \quad m \in M, g \in G.
$$
 For each $m\in M$, let $G(m)= \{\beta_{g}(m): g\in G\}$ be the $G$-orbit of $m$ and $\ds \Omega$ be the {orbit space}, that is, $\Omega = M/ \sim$ with $m\sim n \Leftrightarrow G(m)=G(n)$. Let also:
\begin{itemize}
\item $J_{m}$ be the closed two-sided ideal of $\cA$ generated by maximal ideal $m \in M$ of  $\mathcal{Z} \subset \mathcal{A}$;
\item $\cA_{m}:= \cA / J_{m}$, and $\rho_{m}: \cA \to \cA/J_{m}$ be the quotient map.  
Then by (A1), $J_{\beta_g(m)}= \alpha_g^{-1} J_m$, so it follows that
 $$\cA_{m'}\cong \cA_{m}, \qquad m'\in G(m).$$
\end{itemize}
We  define now a  {family of representations of $\cA$}. For each orbit $\omega\in \Omega$, choose a representative    $m_\omega\in \omega$ and a faithful representation  $\phi_{\omega}: \cA_{m_\omega}\to B(H_{\omega})$ on some Hilbert space $ H_{\omega}$, and define
 $$\pi_{\omega}':  \cA\to B(H_{\omega}), \qquad \pi_{\omega}' = \phi_{\omega} \circ \rho_{m}$$
Recall that to $\pi'_\omega$ we associate $\tilde \pi_{\omega}' : \cA \to B(\ell^2(G, H_\omega))$ such that $(\tilde \pi_{\omega}'(a) \xi)(g):=  \pi_{\omega}' (\alpha_g^{-1}(a))(\xi(g))$. 

\begin{definition} The \emph{local trajectories family on $\cA \rtimes_{\alpha} G$} is the family of regular representations $\{\pi_\omega\}_{\omega \in \Omega}$ induced by $\pi_{\omega}'$, $\omega \in \Omega$, 
that is, induced by the covariant representation $( \pi_{\omega}' ,  \lambda_{\omega})$, such that
$$\pi_{\omega}: \cA\rtimes_{\alpha}G \to B(\ell^{2}(G,H_{\omega})), \qquad  \pi_{\omega}:= \tilde \pi_{\omega}' \rtimes \lambda_{\omega},$$
with
$${\left[\pi_{\omega}(a) \xi\right](t)=\pi_{m_{\omega}}^{\prime}\left(\alpha_{t}^{-1}(a)\right) \xi(t)}, \quad 
{\left[\pi_{\omega}\left(U_{g}\right) \xi\right](t)=\xi\left(g^{-1} t\right)},
$$
for $\xi \in \ell^2(G, H_\omega)$, $t\in G$. Let $\pi= \oplus_{\omega \in \Omega}\pi_{\omega}$ be the direct sum representation of $\cA\rtimes_\alpha G$ on $B(H_\Omega)$, with $H_\Omega= \sum_{\omega\in \Omega} \ell^{2}\left(G, H_{\omega}\right)$. 
\end{definition}

In order for the local trajectories maps to be well-defined on the algebra $\cB=\alg(\cA, U_G)$, we need extra conditions, {namely to guarantee that there is uniqueness of representation of elements in $\cB_0$. 
 One such condition has to do with the structure of fixed points of the action. Recall that we say that  $G$ acts freely on $M$ if the group $\left\{\beta_g: g \in G\right\}$ of homeomorphisms of $M$ onto itself acts freely on $M$, that is, if $\beta_g(m) \neq m$ for all $g \in G \backslash\{e\}$ and all $m \in M$. One considers here a more general notion of freeness that relies on the topology of the state space. 

Let $\cP_{\mathcal{A}}$ be the set of all pure states on $\mathcal{A}$, equipped with the weak-* topology.    As we have noted in \eqref{eq.def.PAtoM}, 
 we have a map 
  \begin{equation}\label{eq.psi}
 \psi:  \cP_\cA\to \hat \cZ\cong M, \quad \mu \mapsto \mu_{|\cZ}.
  \end{equation} 
We often write $m=m_\mu$ if $m=\mu_{|\cZ}= \mathcal{Z} \cap\ker \mu $.

 If $J_m\subset \cA$ is the ideal generated by the maximal ideal $m \in M$ of $\cZ$ and $m=\mu_{|\cZ}$, for some $\mu\in \cP_\cA$,  then $\operatorname{Ker} \mu \supset J_m$ and therefore (see \cite{BKarS_8}, Lemma 4.1) 
\begin{equation}\label{pure}
\cP_{\mathcal{A}}=\bigcup_{m \in M} \mathcal{P}_m, \quad \mathcal{P}_m:= \left\{\mu \in \cP_{\mathcal{A}}:\ker \mu \supset J_m\right\} .
\end{equation}
Writing $m=\mu_{|\cZ}$, then  $\cP_m\subset \cP_\cA$ is the class of extensions of the pure state $\mu_{|\cZ}$ to a pure state on $\cA$.

As in  \cite{Karlovich_OTAA17}, we adopt the following notion of \emph{topologically free action}:

\textit{Condition (A3): For every finite set $G_0 \subset G$ and every open set $W \subset \cP_{\mathcal{A}}$ there exists $\nu \in W$ such that $\beta_g\left(m_\nu\right) \neq m_\nu$ for all $g \in G_0 \backslash\{e\}$, where $m_\nu:=\nu_{|\cZ}=\mathcal{Z} \cap\ker \nu \in M$.
}

 Note that (A3) guarantees that the set of fixed points has empty interior. 
  If the $C^*$-algebra $\mathcal{A}$ is {commutative}, then 
 taking $\mathcal{Z}=\mathcal{A}\cong C(M)$,  we have $\cP_\cA\cong M$, so 
  we can rewrite (A3):

 \emph{(Commutative (A3)) For every finite set $G_0 \subset G$ and every open set $V \subset M$ there exists  $m_0 \in V$ such that $\beta_g\left(m_0\right) \neq m_0$ for all $g \in G_0 \backslash\{e\}$}.

Under assumptions (A1)-(A3), we have an isomorphism $\B\cong  \mathcal{A} \rtimes_\alpha G $ (\cite{Karlovich_OTAA17} - Theorem 3.2). In particular, the local trajectories family is well-defined in $\cB$. Moreover, we have: 
\begin{theorem}[Local trajectories method {- \cite{Karlovich_OTAA17} Theorem 4.1}]\label{thm.mainlocaltrajec}
If (A1)-(A3) hold, 
 then the local trajectories representation $\pi=\oplus_{\omega\in \Omega} \pi_{\omega}$ is faithful in $\B$. Hence, $b\in \B$ is invertible if, and only if, $\pi_{\omega}(b)$ is invertible in $B(\ell^{2}(G,H_{\omega}))$, for all $\omega \in \Omega$, and 
$$
\sup \left\{\left\|\left(\pi_{\omega}(b)\right)^{-1}\right\|: \omega \in \Omega\right\}<\infty.
$$
\end{theorem}

If the number of orbits is {finite}, then the bound on the norms of the inverse elements always holds, so that $b\in \B$ is invertible if, and only if, $\pi_{\omega}(b)$ is invertible in $B(\ell^{2}(G,H_{\omega}))$,} for all $\omega \in \Omega$. In this case, local trajectories family is said to be \emph{sufficient}.

A crucial step in the proof of the above criterion is that the local trajectory family is always injective over $\cA$, that is, $\pi(\cA)\cong \cA$. 
This result relies on the structure of {pure states} of $\cA$ as in \eqref{pure}, we give here a short proof for completeness
 (see also \cite{Karlovich_OTAA17}  for a direct proof of the equality of norms).

\begin{proposition}\label{prop.injectiveA}
Assume (A1) is verified. Then $\{(\pi_{\omega})_{|A}\}_{\omega\in \Omega}$ is a faithful family in $\cA$, i.e., $\pi_{|A}:= \oplus (\pi_{\omega})_{|A}$ is injective and 
$$\|\pi(a)\|=\sup_{\omega \in \Omega} \|\pi_\omega(a)\| = \|a\|_\cA, \quad \mbox{for all } a\in\cA.$$
\end{proposition}
\begin{proof}
Let $\pi(a)= 0$, $a\in \cA$, so that $\pi_\omega(a)=0$, for all $\omega \in \Omega$. 
We show that 
\begin{equation*}\label{eq.kerpiomegaA}
\ker \pi_\omega \cap \cA = \bigcap_{g\in G} J_{\beta_g(m)}, \quad m\in \omega.
\end{equation*}
We have $\pi_\omega(a)= 0 \Leftrightarrow \pi_{\omega}'(\alpha_{g}(a))=0 $, for any $g\in G$ $\Leftrightarrow \rho_{m}(\alpha_g(a)))=0$,  for $g\in G$ $\Leftrightarrow \alpha_g(a) \in J_m$,  for $g\in G$. 
Since, by (A1), $J_m=\alpha_g\left(J_{\beta_g(m)}\right)$, the equality above holds. 

Hence, $\pi(a)= 0$ yields that for any $m\in M$ and $\mu\in \cP_m$, we have $\mu(a)=0$, since $J_m\subset \ker \mu$. By  \eqref{pure}, we have $\mu(a)=0$, for any $\mu\in \cP_\cA$, so $a=0$. 
\end{proof}

Assume (A1) and suppose that $\Phi:= i d \rtimes U: \mathcal{A} \rtimes_{\alpha} G \rightarrow \mathcal{B}$ is an isomorphism,  so the representations $\pi_{w}: \mathcal{B} \rightarrow B\left(\ell^{2}\left(G, H_{w}\right)\right)$ defined by:
$$
{\left[\pi_{\omega}(a) \xi\right](t)=\pi_{m_{\omega}}^{\prime}\left(\alpha_{t}^{-1}(a)\right) \xi(t)}\,, \qquad 
{\left[\pi_{\omega}\left(U_{g}\right) \xi\right](t)=\xi\left(g^{-1} t\right)}
$$
for $\xi \in \ell^{2}\left(G, H_{\omega}\right)$ are well defined. As before, let $\pi= \sum_{\omega\in \Omega} \pi_\omega$.

Now since, by Proposition \ref{prop.injectiveA}, $\pi: \mathcal{A} \rightarrow \pi(\mathcal{A})$ is an isomorphism, we can define $\Psi: C_c(G, \mathcal{A}) \rightarrow$ $C_c(G, \pi(\mathcal{A}))$ given by
$$
\Psi(f)(s)=\pi(f(s)), \quad f\in C_c(G, \mathcal{A}),
$$
which extends to an isomorphism $\Psi: \mathcal{A} \rtimes_\alpha G \rightarrow \pi(\mathcal{A}) \rtimes_{\alpha^{\prime}} G$,  with 
$$\alpha'_g(\pi(a)):= \pi(\alpha_g(a)) = \pi(U_g aU_g^*) = \pi(U_g )\pi(a)\pi(U_g^*). $$
Thus, by definition, $(id, \pi(U))$ is a covariant representation of  $(\pi(\cA),G, \alpha')$. 
We have also that $\pi(\mathcal{B})=\alg\left(\pi(\mathcal{A}), \pi\left(U_G\right)\right)$. 

We then obtain the following commutative diagram:
\begin{equation}\label{CD_loctraj}
\begin{CD}
\mathcal{A} \rtimes_\alpha G @> i d \rtimes_\alpha U>> \cB\\
@VV\cong V @VV\pi V\\
\pi( \mathcal{A}) \rtimes_{\alpha'} G@>i d \rtimes_{\alpha^{\prime}} \pi(U)>> \pi(\cB).
\end{CD}
\end{equation}

The idea of the proof of Theorem \ref{thm.mainlocaltrajec} in \cite{Karlovich_OTAA17} can be summarized roughly as follows: assuming (A1)-(A3) then $ \mathcal{A} \rtimes_{\alpha} G \cong \mathcal{B}$ and moreover $(\pi(\cA), G, \alpha')$ also satisfies (A1) - (A3), since $(id, \pi(U))$ is a covariant representation of  $(\pi(\cA),G, \alpha')$, $\alpha'$ maps $\pi(\cZ)$ to $\pi(\cZ)$, and the condition (A3) of  topological freeness relies on the sets of pure states of $\cA$ and $\pi(\cA)$, which are isomorphic.  This then implies that $\pi$ is an isomorphism due to the commutative diagram above, since all the other arrows are isomorphisms.

\begin{remark}\label{rmk_isom_loctraj} 
 It follows from the discussion above that, assuming (A1), and therefore knowing that $\pi( \mathcal{A}) \cong \mathcal{A}$, with $\pi=\sum_{\omega \in \Omega} \pi_\omega$, then the conclusion of the local trajectories method, Theorem \ref{thm.mainlocaltrajec}, holds in $\cB$ as long as one can show that
$$
i d \rtimes_\alpha U: \mathcal{A} \rtimes_\alpha G \rightarrow \mathcal{B}, \qquad i d \rtimes_{\alpha^{\prime}} \pi(U): \pi(\mathcal{A}) \rtimes_{\alpha^{\prime}} G \rightarrow \pi(\mathcal{B})$$
are isomorphisms, with $\alpha'= \pi\circ \alpha$. The proof of injectivity of $\pi$ on $\cB$ will then follow by the commutativity of the diagram \eqref{CD_loctraj}. Recall that $\cB= \alg(\cA, U_G)$ and $\pi(\cB)=\alg(\pi(\cA), \pi(U_G))$. 
\end{remark}


\section{Isomorphism with the crossed product }\label{s.cond_isom}

As we have noted, one important step to  establish the local trajectory method, Theorem \ref{thm.mainlocaltrajec}, in the non-local algebra $\B= \alg(\cA, U_{G})$. is to give conditions such that $\cB$ is isomorphic to a crossed product algebra. 
Our goal in this section  is to give  conditions that guarantee that such an isomorphism exists, without assuming topological freeness of the action. 

We assume throughout this section  that $\cA\subset B(H)$, for some Hilbert space $H$, $\alpha$ is an action by automorphisms of $\cA$ and that $U: G\to B(H), g\mapsto U_g$ is a unitary representation satisfying
 $$\alpha_{g}(a)= U_{g}a U^{*}_{g},$$ for $g\in G$, $a\in \cA$, that is, that $(id, U)$ is a covariant representation of the dynamical system $(\cA, G, \alpha)$ on $B(H)$. As in Definition \ref{def.algB}, we let $\B:= \alg(\cA, U_{G})$ be the $C^*$-algebra generated by $\cA$ and $U$, which coincides with the closure in $B(H)$ of the $*$-subalgebra
 \begin{equation*}
 \B_0= \left\{ \sum_{g\in G_0} a_g U_g: a_g\in \cA, G_0\subset G \mbox{ finite }\right\}.
 \end{equation*}
The condition that the group $G$ is amenable will be replaced by the weaker condition that the action $\alpha$ is amenable. Note that in the case we have $\cA\rtimes_\alpha G \cong \cA\rtimes_\alpha^r G$. 
\subsection{General conditions and $M$-localization}\label{ss.gencond} 

We start with noting that
 by the universal property of the crossed product, one  has  a surjection $\Phi: C_c(G,\cA) \to \cB_0$ given by $\Phi(f)= \sum_{s\in G} f(s)U_s$ that is bounded in the universal norm and therefore extends to  $\cA \rtimes_{\alpha} G $, so we obtain a surjective map
$$\Phi:= id\rtimes U : \cA \rtimes_{\alpha} G \to \cB.$$
Hence we have  $\cB\cong  \cA\rtimes_{\alpha}G / \ker \Phi$, so
 in order to establish that $\cB \cong \cA \rtimes_{\alpha} G $, it suffices to give conditions such that $\ker\Phi=0$. As it is known, this condition can be written as the boundedness of a family of '{evaluation} maps' \cite{AntonLebedev, Anton_book, Karlovich_OTAA17, RPalma_tese}.

For $s\in G$, let 
\begin{equation}\label{def.E_s}
E_{s}(f):= f(s), \quad f\in C_{c}(G,\cA).
\end{equation} 
Then  $E_s$ is bounded in the universal norm so that it extends to $E_{s}: \cA \rtimes_{\alpha} G \to \cA$. We can write, by continuity, 
$$\Phi(f)= \sum_{s\in G} E_s(f)U_s, \quad f\in \cA \rtimes_{\alpha} G.$$
We will need the following simple lemma.

\begin{lemma}\label{lem.kerEs}
If  $\cA  \rtimes_{\alpha} G\cong  \cA  \rtimes_{\alpha}^{r} G$ then $\bigcap_{s \in G}\ker E_s=\{0\}$. 
\end{lemma}
\begin{proof}
Let $\pi: \cA \to B(H_\cA)$  be an arbitrary   $*$-representation of $\cA$ and $\tilde \pi \rtimes \lambda : \cA \rtimes_{\alpha} G \to  B(\ell^{2}(G,H_\cA))$ be the induced regular representation. Extending \eqref{def.regrep} by continuity, we can write, for $b\in \cA  \rtimes_{\alpha} G$,
$$[(\tilde \pi \rtimes \lambda)(b) \xi](g)=\sum_{s \in G} \pi\left(\alpha_{g^{-1}}[E_s(b)]\right) \xi\left(s^{-1} g\right)$$
for $\xi \in \ell^2(G, H_\cA)$, $g\in G$. If $b\in \bigcap_{s \in G}\ker E_s$, then $ \pi\left(\alpha_{g^{-1}}[E_s(b)]\right)=0$, for all $s, g \in G$, hence $(\tilde \pi \rtimes \lambda)(b)=0$. Since $\pi$ is arbitrary, it follows that for the reduced norm we have  $\|b\|_r=0$, and by the isomorphism $\cA  \rtimes_{\alpha} G\cong  \cA  \rtimes_{\alpha}^{r} G$, also $b=0$. 
\end{proof}
Assume that also $E'_{s}: \B_{0} \to \cA$ such that $$E'_{s}\left(\sum_{g\in G_{0}} a_{g}U_{g}\right) = a_{s}$$ is well defined, that is, we have uniqueness of representation of elements in $\cB_0$. 
Then we have a bijection between $\cB_0$ and $C_c(G,\cA)$ and $E_s(f)=E'_s\circ  \Phi(f) $, for $f\in C_c(G,\cA)$. 
\begin{proposition}\label{prop.conds.isom.crossed}
Assume that the action $\alpha: G\to \Aut(\cA)$  is amenable, in particular, $\cA\rtimes_{\alpha}G\cong \cA\rtimes_{\alpha}^rG$. Then the following statements are equivalent:
\begin{enumerate}[(i)]
\item $\Phi=id\rtimes U: \cA\rtimes_{\alpha}G \to \cB$ is an isomorphism;
\item $E'_{s}: \B_{0} \to \cA$  is well defined and  bounded, for all $s\in G$;
\item For any finite set $G_{0} \subset G$ and $ a_{g}\in \cA$, $g\in G_{0}$, we have
\begin{equation*}
\|a_{e}\| \leq \left\| \sum_{g \in G_{0}} a_{g} U_{g}\right\|.
\end{equation*}
\end{enumerate}
\end{proposition}

\begin{proof}
If  $\Phi=id\rtimes U$ is an isomorphism on $\cA\rtimes_{\alpha}G$, then $E'_s=E_s\circ  \Phi^{-1}$
 is well defined and bounded. Conversely, if $E'_s$ bounded then $E_s= E'_s\circ \Phi$ extends  to $\cA\rtimes_\alpha G$. If $b\in \ker \Phi$ then also $E_s(b)=0$, for all $s$, and since the action of $G$ is amenable, it follows from Lemma \ref{lem.kerEs} that $b=0$. So (i) $\Leftrightarrow$ (ii), and (ii) clearly gives (iii). Assuming (iii), we have $\sum_{g\in G_{0}} a_{g}U_{g}=0 \Rightarrow a_e=0$ and also $a_s=0$, for $s\in G_0$, so we have uniqueness of representation and $E'_s$ is well defined. 
The boundedess of the maps $E'_{s}$ for any $s\in G$ is equivalent to boundedness for $s=e$, since  
$E'_{s}\left(\sum_{g\in G_{0}} a_{g}U_{g}\right) 
= E'_{e}\left(\sum_{g\in G_{0}} a_{g}U_{gs^{-1}}\right)$.
\end{proof}

This result can be found in \cite{AntonLebedev, Anton_book, Karlovich_OTAA17}, with the assumption that $G$ is amenable. 

\begin{remark}\label{rmk.condexp}
The reduced  $C^*$-algebra $ \cA\rtimes_{\alpha}^rG$ always has a canonical  faithful conditional expectation  to $\cA$, that is, a contractive projection onto $\cA$, given by the extension of $E_e$ as in \eqref{def.E_s} (see for instance \cite{BrownOzawa_book}, Proposition 4.1.9). Faithful here means that $E_e(b^*b)=0$ if, and only if, $b=0$ (or equivalently that the kernel of $E_e$ does not contain any non-zero ideal). 
Assuming $\cA\rtimes_{\alpha}G\cong \cA\rtimes_{\alpha}^rG$, we see that the equivalent conditions (ii) and (iii) just mean that there exists a canonical faithful conditional expectation from $\cB= \alg(\cA, U_G)$ to $\cA$ given by $E'_e$. 

\end{remark}

According to the previous proposition, our aim is to give conditions so that we have:

\emph{Condition (B0):
For any finite set $G_{0} \subset G$ and $ a_{g}\in \cA$, $g\in G_{0}$, we have
\begin{equation}\label{B0}
\|a_{e}\| \leq \left\| \sum_{g \in G_{0}} a_{g} U_{g}\right\|.
\end{equation}
}
When \eqref{B0} holds, we say that $(B_0)$ holds for $b= \sum_{g \in G_{0}} a_{g} U_{g}$.

When $\cA=\C$, condition (B0)  gives conditions such that the group $C^*$-algebra $C^*(G)$ and $\alg(U_G)$, the $C^*$-algebra generated by the unitary representations $U_g$, $g\in G$,  are isomorphic.

\begin{remark}\label{rmk.gencovrepB0}
We are working in a concrete setting where $A\subset B(H)$ for some fixed Hilbert space $H$ and $(id, U)$ is a covariant representation. However, condition (B0) and the result above can be used more generally. For an 'abstract' $C^*$-algebra $\cA$ and action $\alpha: G\to \Aut(\cA)$, suppose that we are given any faithful representation $\pi: \cA\to B(H_\pi)$ on a Hilbert space $H_\pi$ and $U: G\to B(H_\pi)$ is a unitary representation. Let $\B_{\pi,U}$ be the $C^*$-algebra generated by $\pi(\cA)$ and $U_g$, $g\in G$. 

Then  $(\pi, U)$ is a covariant representation of $(\cA, G, \alpha)$ on $B(H_\pi)$ if, and only if, $(id, U)$ is a covariant representation of $(\pi(\cA), G, \alpha')$, with $\alpha':= \pi\circ \alpha \circ \pi^{-1}$. Moreover $\alpha$ is amenable if, and only if, $\alpha'$ is amenable. There is an isomorphism $\cA\rtimes_\alpha G\cong \pi(\cA)\rtimes_{\alpha'} G$ that interwines $\pi\rtimes U$
 and $id\rtimes U$. 
 
 We  conclude from Proposition \ref{prop.conds.isom.crossed}  that if  $(\pi, U)$ is a faithful covariant representation of $(\cA, G, \alpha)$ and $\alpha$ is amenable, then $ \pi\rtimes U :\cA\rtimes_\alpha G\to\B_{\pi,U} $ is an isomorphism if, and only if, 
for any finite set $G_{0} \subset G$ and $ a_{g}\in \cA$, $g\in G_{0}$, we have
\begin{equation*}\label{B0gen}
 \left\| \sum_{g \in G_{0}} \pi(a_{g}) U_{g}\right\|\geq \|\pi(a_{e})\|= \|a_{e}\|.
\end{equation*}
\end{remark}

In order to estimate norms, it is often useful to consider positive elements, by property (iii) in Section \ref{ss.rep_states}. We let $b=\sum_{g \in G_{0}} a_{g} U_{g}$, where $G_{0} \subset G$ is a finite set. 
 In terms of pure states, (B0) becomes:
$$\|a_{e}\|^{2}= \|a_{e}^{*}a_{e}\| = \max_{\mu \in \cP_{\cA}} \mu(a_{e}^{*}a_{e}) \leq \max_{ \nu_{0} \in \cP_{B(H)} } \nu_{0}(b^{*}b)= \|b^{*}b\| = \|b\|^{2},$$
We will write often
\begin{equation}\label{b*b}
 b^{*}b= \tilde{a} +\sum_{s \neq t \in G_{0}} \alpha_{s}^{-1}\left(a_{s}^{*} a_{t}\right) U_{s^{-1}t} , \quad \text{with} \quad  \tilde{a}:= \sum_{s \in G_{0}} \alpha_{s}^{-1}\left(a_{s}^{*} a_{s}\right)\in \cA.
\end{equation}
Then $\tilde a$ is positive and  we can write $\tilde{a}= a_e^*a_e + a'$, with $a_e^*a_e$ and $ a'$ positive elements. Hence, $\mu(a_e^*a_e)\leq \mu(\tilde{a})$, for $\mu\in\cP_\cA$, so that also $\|a_e^*a_e\|\leq \|\tilde{a}\|$. It follows that if (B0) holds for all positive elements, that is, if 
$$\|\tilde{a}\|\leq \|b^* b\|, \quad \mbox{ for all  }  b\in \cB_0,$$
then $\|a_{e}\|^{2}= \|a_e^*a_e\|\leq  \|b^* b\|= \|b\|^{2}$ and (B0) holds in general. We conclude that to guarantee that $\Phi$ is an isomorphism, one needs only check (B0) for positive elements $b$ and $a_e$.

 We now give a localized version of (B0). We let $\cZ$ be a central subalgebra  of $\cA$, as in Section  \ref{s.localtraj}.  From now on, we assume condition (A1) of the previous section, that is, that $\alpha$  acts by automorphisms both of $\cA$ and of $\cZ$ and $(id, U)$ is a covariant representation of $(\cA, G, \alpha)$. We will replace condition (A2) of amenability of the group $G$ by the weaker assumption that the action $\alpha$ is amenable. 
 
 As before, we let $M$  be the maximal ideal space of the commutative $C^*$-algebra $\mathcal{Z}$, and identify $\cZ$ with $C(M)$ through the Gelfand transform. 
 Recall from \eqref{eq.def.PAtoM} that 
 we have a surjective, continuous map $$\psi: \cP_\cA \to \cP_\cZ= \hat{\cZ} \cong M, \qquad\mu \mapsto \mu_{|\cZ},$$  
where $\cP_\cA$, $\hat{\cZ}\cong M$ have the weak-$*$ topology. {We have $\cP_\cA= \bigcup_{m\in M} \cP_m$, where $\cP_m$ are the pure states that restrict to $m$ and for $V\subset M$, we write $\cP_V:= \psi^{-1}(V)$.}

We make the following assumption on the algebra $\cA$ of an $M$-localization property:

\emph{Condition (C): For any open set $W\subset \cP_A$, the set  
\begin{equation}\label{C}
V=\psi(W)=\{m\in M: m=\mu_{|\cZ}, \mu\in W\}
\end{equation}
is open in $M\cong \hat{\cZ}$.
}

We say that such an algebra $\cA$ is \emph{$M$-localizable}. Clearly any commutative algebra satisfies (C), since in this case $\cP_\cA$ is homeomorphic to $M$. Also matrix algebras $\left[C(M)\right]_{N\times N}$, $N\in \N$, satisfy (C), and more generally also algebras of the form $HOM(E,F)$, where $E, F$ are vector bundles, as the ones considered in \cite{Anton_book}. Moreover, all $C^*$-algebras that have the uniqueness of extension property \cite{anderson, archbold_ext_states, ChuKusuda_ext_states}), that is, if $\psi$ is injective in that any pure state in  $\cZ$ has a unique extension to a pure state in $\cA$, also satisfy (C), since in this case the extension in $\cS_\cA$ is also unique and the extension map $\cP_\cZ \to \cP_\cA \subset \cS_\cA$ is continuous (see \cite{ChuKusuda_ext_states}, 
Lemma 1, or directly).

Note that if $\cA$ satisfies (C), then topological freeness, condition (A3), becomes equivalent to requiring that no non-empty open set of $M$ is fixed by a finite, non-trivial, subset $G_0\subset G$, similarly to the commutative case and to the definition adopted in \cite{Anton_book}.

According to \eqref{b*b}, (B0) for the positive element $b^*b$ comes down to $\|b^*b\|\geq \|\tilde a\|$. The next result shows that, assuming our algebra satisfies (C),  if this inequality holds locally, for any positive element, then also (B0)  holds in $\cB$.

We first  introduce the following notation that we shall use throughout the paper: for  $\emptyset \neq V\subset M$ open, let $\cZ(V)\subset \cZ$ be those $z\in C(M)$ with $0\leq z\leq 1$, $\|z\|=1$, $\supp z \subset V$. We think of elements $zb$, $z\in \cZ(V)$, $b\in \cB$, as being a localization of $b$ to $V$. 

 If we let $\mu\in \cP_\cA$, with $m_\mu=\mu_{|\cZ}\in M$ and $z\in \cZ$, then if $z(m_\mu)=0$ then $\mu_0(zb)=0$, for all $b\in \cB$, $\mu_0$ any extension of $\mu$ to a state in $\cB$, since
$$|\mu_0(zb)|^2\leq \mu(zz^*) \mu_0(b^*b)= z^2(m_\mu)\mu_0(b^*b)=0.$$
 In particular, if $V\subset M$ is open and $z_V\in \cZ(V)$, then  for  $\mu_0$ an extension of $\mu$ with $\mu\notin \cP_V$, that is, $m_\mu\notin V$, we have 
 \begin{equation}\label{eq.muzVb=0}
 \mu_0(z_Vb)=0, \quad \mbox{for }b\in \cB.
 \end{equation}
Moreover,  if $a\in \cA$ with $a\geq 0$ then, since $a$ and $z_V\in \cZ$ commute and $z_V\geq 0$, we have $z_Va\geq 0$ and
$$\|z_Va\|= \max_{\mu\in \cP_V} \mu(z_Va)=\max_{\mu\in \cP_V} z_V(\mu_{|\cZ})\mu(a).$$ 
We  now prove:
\begin{proposition}\label{localB1}
Assume $\cA$ satisfies (C). Let $G_0\subset G$ be finite and $b= \sum_{g\in G_0} a_g U_g$, $a_g\in \cA$, and $\tilde a : = \sum_{g\in G_0} \alpha_g^{-1}(a_g^*a_g)$. If for any non-empty open set $V\subset M$, there is a central element $z_{V}\in \cZ(V)$
 satisfying 
\begin{equation}\label{eq.B1+}
\| z_{V} b^{*}b \| \geq \| z_{V} \tilde{a}\|
\end{equation}
then $(B0)$ holds for $b$. 
In fact, it suffices that $\| z_{V} b^{*}b \| \geq \nu( z_{V} \tilde{a})$, for $\nu\in \cP_{\cA}$ such that $z_{V}(m_{\nu})=1$, with $m_\nu=\nu_{|\cZ}$.
\end{proposition}

\begin{proof}
Let $\phi\in \cP_{\cA}$ be such that $\|a_{e}^{*}a_{e}\|= \phi(a_{e}^{*}a_{e})\leq \phi(\tilde{a})$ and consider the open sets 
$$W=\{ \mu\in \cP_\cA: |\mu(\tilde{a}) - \phi(\tilde{a})|<\eps\} \quad \mbox{  and  }
\quad V=\{\mu_{|\cZ}: \mu\in W\} \subset M.$$
Let $z_{V}\in \cZ(V)$. Then since $z_{V} \tilde{a}$ is positive, 
writing $m_\mu=\mu_{|\cZ}$, $\mu\in \cP_\cA$, 
\begin{align*}
 \| z_{V} \tilde{a}\| & = \max_{\mu\in \cP_{\cA}} \mu(z_{V})\mu(\tilde a)
\geq \max_{\mu\in W} \mu(z_{V})\mu(\tilde a)\\
& >  ( \phi(\tilde{a}) - \eps) \max_{m_\mu\in V} z_{V}(m_\mu)= \phi(\tilde{a})-\eps\geq \|a_{e}\|^{2}-\eps.
\end{align*}
Taking $z_V$ satisfying the assumption,  and noting that $\| z_{V} b^{*}b \| \leq \|b\|^{2}$, we get  $\|b\|^{2}\geq\|a_{e}\|^{2}-\eps $, for all $\eps>0$, so (B0) follows. 

For the last assertion, in this case, taking $m\in V$ with $z(m)=1$ and $\nu\in W$ with $\nu_{|\cZ}= m$, we get  $\| z_{V} b^{*}b \| \geq \nu( \tilde{a})$ for some $\nu\in W$,  and proceed in the same way.
\end{proof}

 In terms of pure states,  the condition in Proposition \ref{localB1} follows if for any open $V\subset M$, there is $z_{V}\in \cZ(V)$  and  {$\nu\in \cP_\cA$ such that}, 
$$\nu_{0}\left( z_{V}\tilde{a} + \sum_{s \neq t \in G_{0}} \alpha_{s}^{-1}\left(a_{s}^{*} a_{t}\right)  z_{V}U_{s^{-1}t}\right) \geq \mu(z_{V} \tilde{a}), \quad \mbox{  for  all  } \mu \in \cP_{\cA}, {\mbox{  with  } \mu_{|\cZ}\in V},$$
and $\nu_{0}$ an extension of $\nu$ to a state in $\cB$. In fact it suffices to check this for $\mu\in\cP_{\cA}$ with $\mu(z_{V})=1$. 

We arrive at condition (B1):

\emph{Condition (B1): For any non-empty open set $V\subset M$, finite set $G_0\subset G$ and $a_g \in \cA$, $g\in G_0$, there exists a central element $z_{V}\in \cZ(V)$
 satisfying 
\begin{equation}\label{eq.B1}
 \| z_{V} {a_e}\| \leq \left \| z_{V} \sum_{g\in G_0} a_g U_g \right \|.
\end{equation}
}

As a matter of terminology, if \eqref{eq.B1} or \eqref{eq.B1+} hold for $b= \sum_{g\in G_0} a_g U_g\in \cB_0$, we say that $(B1)$ holds for $b$, and if it holds on a subset $\cB_0'$ dense in a $C^*$-subalgebra $\cB'\subset \cB$ we say that (B1) holds in $\cB'$. 

Of course, if (B0) holds, then (B1) holds {for any} $z_V\in \cZ(V)$. On the other hand, the condition in Proposition \ref{localB1} amounts to (B1) for the positive elements $b^*b$ and $\tilde{a}$. We have then that for $M$-localizable algebras $\cA$, as in \eqref{C}, the global condition (B0) is equivalent to the $M$-local condition (B1). 

Together with Proposition \ref{prop.conds.isom.crossed}, we obtain another condition for the isomorphism of $\cB$ with the crossed product algebra. 
 Recall that condition (A1) means  that $(id, U)$ is a covariant representation of $(\cA, G, \alpha)$, where $\alpha$  acts by automorphisms both of $\cA$ and of the central subalgebra $\cZ$.  

\begin{corollary}\label{cor.B1thenB0}
Assume that $\cA$ verifies (C), condition (A1) is verified for $(\cA, G, \alpha)$ and that the action $\alpha$ is amenable. Then
$$ (B1)\quad  \Leftrightarrow \quad (B0)\quad  \Leftrightarrow \quad \Phi=id\rtimes U : \cA\rtimes_\alpha G \to \cB \mbox{  is an isomorphism}.$$
\end{corollary}

In the spirit of Remark \ref{rmk.gencovrepB0}, if we are given an arbitrary faithful representation $\pi: \cA\to B(H_\pi)$, with $H_\pi$ a Hilbert space, and a covariant representation $(\pi, U)$, with $U: G \to B(H_\pi)$ a unitary representation, we  can also easily deduce a $M$-local condition to ensure that $\cA \rtimes_\alpha G\cong \cB_{\pi, U}:= \alg(\pi (\cA), U_G)$. 

We will see in the next section that, in fact, for these results to hold, we only need to check (B1) for a subclass of open subsets of $M$.

\begin{remark}\label{rmk.isomdeprep}
We remark that conditions (B0) and (B1) are not, in general, related to the action $\alpha$ but only to the way $\cA$ interacts with the unitaries $U_g$. In some situations, we may have unitaries $U_g$ and $U_g'$ defining the same action $\alpha :G\to \Aut(\cA)$ such that the algebras $\cB= \alg(\cA, U_G)$ and $\cB'= \alg(\cA, U_G')$ are {not} isomorphic, if (B0) is satisfied in $\cB$ but not in $\cB'$ (see for instance Example 12.11 in \cite{AntonLebedev}). This cannot happen if $G$ is amenable and the action is \emph{topologically free}, by the isomorphism theorems in \cite{Karlovich_OTAA17} (Theorem 3.3) and \cite{AntonLebedev} (Corollary 12.16), as in this case we would have $\cB\cong \cB'\cong \cA\rtimes_\alpha G$.

We note further that if the action is \emph{not} topologically free, then the fact that conditions (B0) and (B1) hold may depend on the way $\cA$ is represented as an algebra of bounded operators;
 see Examples \ref{example1} and \ref{example2} at the end of the next section.
\end{remark}


\subsection{Fixed points}\label{ss.fixed_points}

We show here how the structure of fixed points plays a role in ensuring the localized condition (B1).  We keep the assumption  that  condition (A1) holds for our dynamical system $(\cA, G, \alpha)$, in that the action $\alpha$  leaves the central subalgebra $\cZ= C(M)$ invariant and $(id, U)$ is a covariant representation of $(\cA, G, \alpha)$ on some Hilbert space $H$, with  $U$ a unitary representation of $G$ in $B(H)$, and $\cB= \alg(\cA, U_G)$. We also consider the induced action of $G$ on $M$ given by 
$ \beta_{g}: M\to M $ %
 such that $
z\left(\beta_g(m)\right)=\left(\alpha_g(z)\right)(m)$, $z \in \mathcal{Z}, m \in M, g \in G.
$

 Note that given an open set $V\subset M$ and finite $G_0\subset G$,  if $m_0\in V$  is such that $\beta_{g}(m_0)\neq m_0$, for all $g\in G_{0}$, then there exists an open $\Delta\subset V$, $m_0\in \Delta$, such that  
 \begin{equation}\label{eq.Delta}
 \beta_{g}(m)\neq m, \quad \text{for all} \quad m\in \Delta, g\in G_{0}, \quad \text{and} \quad \beta_{g}(\Delta)\cap \Delta = \emptyset.
 \end{equation}
In general, we have the the following lemma.
\begin{lemma}\label{lem_fixedpts}
 Given a finite set $G_{0} \subset G$ and a non-empty open set $V \subset M$, there exists a non-empty open set $\Delta \subset V$ such that for $g\in G_{0}$, either 
  $\left.\beta_{g}\right|_{\Delta} = i d_{\Delta}$ or  $\beta_{g}(\Delta) \cap \Delta=\emptyset$.
 In particular, we have $G_{0}= \tilde{D} \cup D_{0}$ such that 
$$\left.\beta_{g}\right|_{\Delta}=i d_{\Delta}, \quad\mbox{ for }  g \in \tilde{D},\qquad \beta_{g}(\Delta) \cap \Delta=\emptyset, \quad  \mbox{ for } g\in D_{0}=G_{0}-\tilde{D}.$$
 Moreover,{ if $G_0$ is closed for inverses}, then $D_{0}$ and $\tilde{D}$ are closed for inverses. 
\end{lemma}

\begin{proof} Induction on number of elements of $G_{0}$: if $G_0=\{g\}$, then either  $\left.\beta_{g}\right|_{V} = i d_{V}$, and $\tilde{D}:= G_0$, or there exists $m\in V$ with $\beta_g(m)\neq m$,  in which case we can take $\Delta\subset V$ open with $\beta_{g}(\Delta) \cap \Delta=\emptyset$ and $D_0=G_0$.

 Assuming now that $G_0=\{g_k : 0\leq k \leq n\}$, let $\Delta' \subset V$ be such that $G_{0} \setminus \{g_n\}= \tilde{D}' \cup D'_{0}$ with
$\left.\beta_{g}\right|_{\Delta'}=i d_{\Delta'},$ for $  g \in \tilde{D}'$ and $ \beta_{g}(\Delta') \cap \Delta'=\emptyset,$ for $ g\in D_{0}'.$ As above, if
 $\left.\beta_{g_n}\right|_{\Delta'} = i d_{\Delta'}$, we take $\Delta=\Delta'$, $\tilde{D}= \tilde{D}' \cup\{g_n\}$, $D_0=D_0'$. Otherwise, there exists $m\in \Delta'$ with $\beta_{g_n}(m)\neq m$,  in which case we can take $\Delta\subset \Delta'$ open with $\beta_{g_n}(\Delta) \cap \Delta=\emptyset$ and $D_0 = D_0' \cup\{g_n\}$.

Assuming $G_0$ closed for inverses, then   $\left.\beta_{g}\right|_{\Delta} = i d_{\Delta} \Leftrightarrow \left.\beta_{g^{-1}}\right|_{\Delta} = i d_{\Delta} $ and  $\beta_{g}(\Delta) \cap \Delta=\emptyset \Leftrightarrow \beta_{g^{-1}}(\Delta) \cap \Delta =\emptyset$, so also $D_0$ and $\tilde{D}$ are closed for inverses.
\end{proof}

Our goal now is to see that to prove condition (B1), it suffices to consider open sets of $M$ that are fixed by a finite subset of $G$. 

\begin{remark}\label{rmk.opensets}
For $z\in \cZ$ and $g\in G$, we have $zU_g = U_g \alpha_g^{-1}(z) = U_g (z\circ \beta_{g^{-1}})$. Let $\Delta\subset M$ be open and $z_\Delta\in \cZ(\Delta)\subset C(M)$, that is, $z_\Delta$ is a non-negative function, supported in $\Delta$, with $\|z_\Delta\|=1.$ We will use the following.
\begin{enumerate}[(i)]
\item If $\left.\beta_{g}\right|_{\Delta} = i d_{\Delta}$, then  $z_\Delta U_g =U_gz_\Delta$, that is, $z_\Delta$ and $U_g$ commute. (This holds as long as $z_\Delta \in \cZ$ is zero outside $\Delta$.)

\item If  $\beta_{g}(\Delta) \cap \Delta=\emptyset$, then $z_\Delta \circ \beta_{g^{-1}}$ is zero on $\Delta$.
 In particular, if $\mu_0$ is a state of $\cB$ that restricts to a pure state in $\cA$, so that ${\mu_0}_{|\cZ}= m_\mu\in \cZ$, then $\mu_0(z_\Delta aU_g)=0$: if $m_\mu \notin \Delta$ it follows from \eqref{eq.muzVb=0}, if $m_\mu\in \Delta$ then
$$|\mu_0(z_\Delta aU_g)|^2= |\mu_0( aU_g (z_\Delta \circ \beta_{g^{-1}}))|^2 \leq \mu(aU_gU_g^*a^*)(z_\Delta \circ \beta_{g^{-1}})^2(m_\mu)=0.$$
\end{enumerate}
\end{remark}

\begin{proposition}\label{prop.B2prev}
Given $G_{0}\subset G$ finite and a non-empty open set  $V\subset M$, let $G_0'=\{s^{-1}t: s,t\in G_0\}= D_0 \cup \tilde{D},$  and $\emptyset\neq \Delta \subset V$ be an open set
 such that $\left.\beta_{g}\right|_{\Delta}=\left.i d\right|_{\Delta}$ for $g=s^{-1}t\in \tilde D$, and $ \beta_{g}(\Delta) \cap \Delta=\emptyset, $ for $g=s^{-1}t \in D_{0}=G'_{0}-\tilde{D}$, $s,t  \in G_{0}$. 
 
 Then for  $b= \sum_{g\in G_{0}} a_{g}U_{g}$, $a_g\in \cA$,  and any $z_\Delta\in \cZ(\Delta)$, we have that
$$\|z_\Delta b^* b \|\geq \left\|z_{\Delta}\left( \tilde{a} + \sum_{s \neq t, s^{-1}t \in \tilde{D}} \alpha_{s}^{-1}\left(a_{s}^{*} a_{t}\right) U_{s^{-1}t}\right) \right\|,$$
with $ \tilde{a}=  \sum_{s \in {G_0}} \alpha_{s}^{-1}\left(a_{s}^{*} a_{s}\right) $.
\end{proposition}

Note that the existence of such an open set $\Delta$ and sets $ D_0, \tilde{D}\subset G$ is a consequence of applying Lemma \ref{lem_fixedpts} to the finite set $G_0'=\{s^{-1}t: s,t\in G_0\}.$
\begin{proof}
Let $\emptyset\neq V\subset M$ be  open and  let $b= \sum_{g \in G_{0}} a_{g} U_{g}$, where we assume without loss of generality that $G_0$ is closed for inverses, so also $G_0'$ is closed for inverses. Write $G_0'= D_0 \cup \tilde{D}$ satisfying the conditions above, with $\Delta \subset M$ open. 

 Then we can write
 \begin{align*}
 b^{*}b & =\sum_{s , t \in G_{0}} \alpha_{s}^{-1}\left(a_{s}^{*} a_{t}\right) U_{s^{-1}t}\\
 & = \tilde{a} + \sum_{s \neq t,  s^{-1}t\in {D_0}} \alpha_{s}^{-1}\left(a_{s}^{*} a_{t}\right) U_{s^{-1}t}  +\sum_{s \neq t,  s^{-1}t\in \tilde{D}} \alpha_{s}^{-1}\left(a_{s}^{*} a_{t}\right) U_{s^{-1}t}.
 \end{align*}
  Let $\tilde{b}:= \tilde{a} + \sum_{s \neq t, s^{-1}t \in \tilde{D}} \alpha_{s}^{-1}\left(a_{s}^{*} a_{t}\right) U_{s^{-1}t}$. 
  Note that $\tilde{b}$  is self-adjoint since $\tilde{D}$ is closed for inverses, and 
$$\left( \alpha_{s}^{-1}\left(a_{s}^{*} a_{t}\right) U_{s^{-1}t}\right)^*= \alpha_{t}^{-1}\left(a_{t}^{*} a_{s}\right) U_{t^{-1}s}.
$$
Consider  the $C^*$-subalgebra ${\mathcal{Z}_0(\Delta)} \subset \mathcal{Z}=C(M)$ of  functions that are zero outside $\Delta$ and let $z_\Delta \in \cZ_0(\Delta)$. 
Since $\left.\beta_{g}\right|_{\Delta}=i d_{\Delta},$ for  $g= s^{-1}t \in \tilde{D}$, we have that $z_\Delta U_g=U_g (z_\Delta \circ \beta_{g^{-1}})= U_g z_\Delta$, from which follows that  $z_\Delta$ and $\tilde b$ commute.

In particular,  $z_\Delta \tilde b$ is also self-adjoint, hence normal, so there always exists a pure state $\nu_\cB\in \cP_\cB$  such that
 $\|z_\Delta \tilde b\|= |\nu_\B(z_\Delta \tilde b)|.$  We see now that we can pick such a $\nu_\B$ such that $\nu_\cB(z_\Delta b^*b)= \nu_\cB(z_\Delta \tilde b)$, which proves our claim since in this case
 $$
\left\|z_{\Delta} b^* b\right\| \geq | \nu_\cB(z_{\Delta} b^*b)| = | \nu_\cB(z_{\Delta}\tilde b)|= \|z_\Delta \tilde b\|.
$$
 
Consider then the $C^*$-subalgebras
$$
\mathcal{C}:=\operatorname{alg}\{\tilde b,  {\mathcal{Z}_0(\Delta)}, I d\} \subset \mathcal{B} \qquad \mbox{and} \qquad \tilde{\mathcal{C}}: =\operatorname{alg}\{\mathcal{C}, \mathcal{Z}\}.
$$
Let $z_{\Delta} \in \mathcal{Z}(\Delta)$, that is, $z_{\Delta} \in {\mathcal{Z}_0(\Delta)} $ such that $\left\|z_{\Delta}\right\|=1$, and assume $z_\Delta \tilde b\neq 0$ (otherwise there is nothing to prove). 
 Since $\cC$ is a commutative $C^*$-algebra, as $\tilde b$ commutes with any element of $Z_0(\Delta)$, there exists a pure state $\nu$ of $\cC$
 such that $$\|z_\Delta \tilde b\|= |\nu(z_\Delta \tilde b)|.$$
Moreover, $\nu(z_\Delta \tilde b)= \nu(z_\Delta) \nu( \tilde b)\neq 0$, so $\nu(z_\Delta)\neq 0$.

Since $\mathcal{C} \subset \tilde{\mathcal{C}}$ is a subalgebra and $\nu$ is a pure state of $\mathcal{C}$, there is an extension of $\nu$ that is a pure state of $\tilde{\mathcal{C}}$, $\nu_{\tilde{\mathcal{C}}} \in \mathcal{P}_{\tilde{\mathcal{C}}}$. Let $f \in \mathcal{Z}$ such that $f z_{\Delta}=0$. Then since $z_{\Delta} \in \tilde{\mathcal{C}}$ is a central element of this subalgebra, then
$$
0=\nu_{\tilde{C}}\left(z_{\Delta} f\right)=\nu_{\tilde{C}}\left(z_{\Delta}\right) \nu_{\tilde{C}}(f)
$$
and since $\nu_{\tilde{C}}\left(z_{\Delta}\right)=\nu\left(z_{\Delta}\right) \neq 0$ then $\nu_{\tilde{C}}(f)=0$. Thus we have that
$$
\nu_{\tilde{C}}(f)=0 \quad \mbox{  for all  }f \in \mathcal{Z}: f z_{\Delta}=0
$$
Since $\tilde{C} \subset \mathcal{B}$ is a subalgebra, we can extend the pure state $\nu_{\tilde{C}}$ to a pure state $\nu_{\mathcal{B}} \in \mathcal{P}_{\mathcal{B}}$. Then if $s^{-1} t \in D_0$ we have that
$$
z_{\Delta} \alpha_s^{-1}\left(a_s^* a_t\right) U_{s^{-1} t}=\alpha_s^{-1}\left(a_s^* a_t\right) U_{s^{-1} t}\left(z_\Delta \circ \beta_{t^{-1} s}\right)
$$
and since  $ \beta_{g}(\Delta) \cap \Delta=\emptyset, $ for $g=s^{-1}t \in D_{0}$, we have that
$$
\left(z_\Delta \circ \beta_{t^{-1} s}\right)^2 z_{\Delta}=0
$$
and so
$$
\nu_{\mathcal{B}}\left(\left(z_\Delta \circ \beta_{t^{-1} s}\right)^2\right)=\nu_{\tilde{\mathcal{C}}}\left(\left(z_\Delta \circ \beta_{t^{-1} s}\right)^2\right)=0
$$
for $s \neq t$ and $s^{-1}t \in D_{0}$. 
By Cauchy-Schwartz inequality for states, we have that:
$$
\begin{aligned}
\left|\nu_{\mathcal{B}}\left(z_{\Delta} \alpha_s^{-1}\left(a_s^* a_t\right) U_{s^{-1} t}\right)\right|^2 & =\left|\nu_{\mathcal{B}}\left(\alpha_s^{-1}\left(a_s^* a_t\right) U_{s^{-1} t}\left(z_\Delta \circ \beta_{t^{-1} s}\right)\right)\right|^2 \\
& \leq \nu_{\mathcal{B}}\left(\alpha_s^{-1}\left(a_s^* a_t\right) U_{s^{-1} t}\left(\alpha_s^{-1}\left(a_s^* a_t\right) U_{s^{-1} t}\right)^*\right) \nu_{\mathcal{B}}\left(\left(z_\Delta \circ \beta_{t^{-1} s}\right)^2\right)\\
&=0 .
\end{aligned}
$$
Thus we have that
$$
\left|\nu_{\mathcal{B}}\left(z_{\Delta} b^* b\right)\right|=\left|\nu\left(z_{\Delta}\tilde b\right)\right|=\| z_{\Delta}\tilde b\|
$$
and from this we conclude that
$
\left\|z_{\Delta} b^* b\right\| \geq \| z_{\Delta}\tilde b\|.
$
\end{proof}

One first consequence applies when the action is topologically free, as in condition (A3) of Section \ref{s.localtraj}. It follows from \eqref{eq.Delta} that if  (A3) holds, that is, if
 the action is topologically free, then for any open $V\subset M$, and $G_0\subset G$ finite, we can always find  $\Delta\subset V$ with $\beta_{g}(\Delta)\cap \Delta = \emptyset$, for all $g\in G_0$.
 Hence, in Lemma  \ref{lem_fixedpts} we can always take $D_0=G_0$ and $\tilde D=\emptyset$.

We obtain the following version of the result  in \cite{Karlovich_OTAA17} (Theorem 3.2) that establishes the isomorphism of $\cB$ with the crossed product . Recall that we always assume condition (A1), so that in particular $(id, U)$ is a covariant representation of $(A, G, \alpha)$, and that (C) holds if $\cA$ is $M$-localizable as in \eqref{C}. 

\begin{corollary}
Assume  that (A1) is verified and that the action is topologically free, in that (A3) is verified, then (B1) holds in $\cB$. 
If moreover (C) holds for $\cA$ and the action $\alpha$ is amenable,  then $id\rtimes U : A\rtimes_\alpha G \to \cB$ is an isomorphism.
\end{corollary}

\begin{proof}
Let  $b= \sum_{g\in G_{0}} a_{g}U_{g}$, $a_g\in \cA$, $G_0$ finite, and $\emptyset \neq V\subset M$ be open. If $\cA$ satisfies (C) and (A3) holds, we can take $\tilde D= \emptyset$ in Proposition \ref{prop.B2prev}, so we get for any $z_\Delta\in \cZ(\Delta)$, 
$$\|z_\Delta b^* b \|\geq \left\|z_{\Delta}\tilde{a}  \right\|,$$
with $ \tilde{a}=  \sum_{s \in {G_0}} \alpha_{s}^{-1}\left(a_{s}^{*} a_{s}\right) $.
We conclude that (B1) holds for $b$. It follows from Corollary \ref{cor.B1thenB0} that  if the action $\alpha$ is amenable, then $id\rtimes U$ is an isomorphism  .
\end{proof}

Even if the action is not topologically free, we can still get a criterion for the isomorphism of $\cB$ with the crossed product to hold. The following result comes directly from Proposition \ref{prop.B2prev}, as in this case, \eqref{eq.B1+}  holds. 

\begin{corollary}\label{prop.B2}
Let $G_{0}\subset G$ be finite, $a_g\in \cA$, $g\in G_0$, and $\tilde a  = \sum_{g\in G_0} \alpha_g^{-1}(a_g^*a_g)$. 
If for every non-empty open set $V\subset M$ such that $\left.\beta_{g}\right|_{V}=\left.i d\right|_{V}$ for $g\in D$, with ${D}\subset \{s^{-1}t:  s,t  \in G_{0}\}$ arbitrary, there is a central element $z_{V}\in \cZ(V)$  
 satisfying 
$$\left\|z_{V}\left( \tilde{a} + \sum_{s \neq t, s^{-1}t \in {D}} \alpha_{s}^{-1}\left(a_{s}^{*} a_{t}\right) U_{s^{-1}t}\right) \right\|\geq \| z_{V} \tilde{a}\|$$
then $(B1)$ holds for $b= \sum_{G_{0}} a_{g}U_{g}$. 
\end{corollary}

We arrive at condition (B2):

\emph{Condition (B2):  For any finite set $D \subset G$ and any non-empty open $V \subset M$ such that $\left.\beta_{g}\right|_{V}=\left.i d\right|_{V}$ for all $g \in D$, and  $a_{g} \in \cA$, $g\in D$, there exists a central element  $z_{V} \in \mathcal{Z}(V)$  
 satisfying 
 \begin{equation}\label{eq.B2}
\left\|z_{V}{a_e}\right\| \leq  \left\|z_{V}\sum_{g  \in D} a_{g} U_{g}\right\|.
\end{equation}
}

We can assume that ${a_e}$ is positive and $\sum_{g \neq e \in D} a_{g} U_{g}$ is self-adjoint. If \eqref{eq.B2} holds for $b= \sum_{g\in G_0} a_g U_g\in \cB_0$, we say that (B2) holds for $b$, and if it holds on a subset $\cB_0'$ dense in a $C^*$-subalgebra $\cB'\subset \cB$ we say that (B2) holds in $\cB'$. 

In the notation of Corollary \ref{prop.B2},  applying (B2) to the element $$\tilde{b}= \tilde{a} + \sum_{s \neq t, s^{-1}t \in {D}} \alpha_{s}^{-1}\left(a_{s}^{*} a_{t}\right) U_{s^{-1}t},$$ we obtain that (B1) holds for $b= \sum_{G_{0}} a_{g}U_{g}.$
Of course, if (B1) holds then we can apply it 
 on open sets $V$ such that $\left.\beta_{g}\right|_{V}=\left.i d\right|_{V}$ for all $g \in D$, to obtain (B2). Note that such sets exist only when the action is not topologically free. 

 It follows from the previous corollary that we can improve Corollary \ref{cor.B1thenB0} to get conditions for the isomorphism of $\cB= \alg(\cA, U_G)$ with the crossed product algebra. 
 Recall that (C) means that $\cA$ is $M$-localizable as in \eqref{C} and that condition (A1) is simply  that $(id, U)$ is a covariant representation of $(\cA, G, \alpha)$, where $\alpha$  leaves the central subalgebra $\cZ$ invariant.

\begin{corollary}\label{cor.B2equivB1}
Assume that $\cA$ verifies (C), that condition (A1) is verified for $(\cA, G, \alpha)$ and that the action $\alpha$ is amenable. Then 
$$(B1) \quad   \Leftrightarrow \quad (B2)\quad  \Leftrightarrow \quad \Phi=id\rtimes U : \cA\rtimes_\alpha G \to \cB \mbox{  is an isomorphism}.$$
\end{corollary}

In the setting of arbitrary faithful representations, as in Remark \ref{rmk.gencovrepB0}, if we are given a faithful representation $\pi: \cA\to B(H_\pi)$, with $H_\pi$ a Hilbert space, and a covariant representation $(\pi, U)$ of $(\cA, G, \alpha)$, with $U: G \to B(H_\pi)$ a unitary representation, we  can replace \eqref{eq.B2} in (B2) by
$$\left\|z_{V}a_e\right\| \leq  \left\|\pi(z_{V})\sum_{g  \in D} \pi(a_{g}) U_{g}\right\|,$$
and, assuming that $\cA$ satisfies (C) and the action $\alpha$ is amenable, we can conclude that $\cA \rtimes_\alpha G\cong \cB_{\pi, U}:= \alg(\pi (\cA), U_G)$.

We now apply our results to show that if topological freeness of the action is not satisfied for some system $(\mathcal{A}, G, \alpha)$, it is possible to construct covariant representations $(\pi,U)$ for which (B2) is not satisfied and so the $C^*$ algebra $\text{alg}(\mathcal{A}, U_G)$ is not isomorphic to the crossed product. (Note that for a commutative algebra, the existence of such a covariant representation follows from Theorem 2 in \cite{ArchboldSpielberg}.)

\begin{example}\label{example1}
   Let $\mathcal{A} = C(M)$ for some compact $M\subset \mathbb{R}^N$ such that $M$ is the closure of an open set and $M$ has smooth boundary,
 and  $G$ be an amenable discrete group. Let
    $\beta : G \times M \rightarrow M$ be an action by diffeomorphisms such that $\beta_e=id$ and $\beta_{gs}= \beta_g\beta_s$, $g, s\in G$. Assume that for some $g \in G$ where $g \neq e$ there exists a non-empty open set $V \subset M$ such that
    $$
        \beta_g|_V = id|_V.
    $$
    Let $\alpha: G\to \Aut(\cA)$ be given by $\alpha_g(a) = a\circ \beta_g$, $g\in G$. 
    We claim that there exists a covariant representation $(\pi,U)$ of $(\mathcal{A}, G, \alpha)$ such that condition (B2) is not satisfied in $\cB= \alg(\pi(\cA), U)$.
    
     Indeed, consider the Hilbert space given by $H = L^2(M)$ with the Lebesgue measure. Consider the faithful covariant representation on $B(H)$ given by
    \begin{equation*}
        [\pi(a)f](x)
    =
        a(x)
        f(x),
    \quad\quad
        [U_gf](x)
    =
        |\det d\beta_g|^{1/2}
        f(\beta_g(x)),
    \end{equation*}
    for $a\in \cA=C(M)$ and $f\in L^2(M)$. 
    Let  $e\neq g\in G$ and $\emptyset \neq V\subset M$ be such that $
        \beta_g|_V = id|_V,$  and let $z_V $ be a continuous  non-negative function on $M$ with compact support in $V$ with $\|z_V\|=1$. Then since $\beta_g|_V = id|_V$, we have
    \begin{align*}
      [  (\pi(z_V)
        U_g) f](x)
    =
        z_V(x)  |\det d\beta_g|^{1/2}
        f(\beta_g(x))
    =
        z_V(x)
        f(x)
    =
        [\pi(z_V)f](x).
    \end{align*}
    Thus
    \begin{equation*}
        0=\|\pi(z_V) - \pi(z_V)U_g\|
        <
        \|\pi(z_V)\|=\|z_V\|=1,
    \end{equation*}
    and so (B2), hence also (B0), is not satisfied. Since $\cA$ is $M$-localizable, as it is commutative, and $G$ is amenable, we conclude from Corollary \ref{cor.B2equivB1}, or directly from Proposition \ref{prop.conds.isom.crossed} (see also Remark \ref{rmk.gencovrepB0}), that $\cB$ is not isomorphic to the crossed product $\cA\rtimes_\alpha G$. Note that this does not contradict the isomorphism theorems in \cite{AntonLebedev, Karlovich_OTAA17} as the action $\alpha$ is not topologically free, 
\end{example}

Although the above example shows that for such an algebra $\cA$ and action $\alpha$ some covariant representations of $(A, G, \alpha)$ may not satisfy condition (B0), however it is possible to construct concrete representations which satisfy it.

\begin{example}\label{example2}
    Consider the $C^*$-algebra $\mathcal{A} = C(M)$  for some compact $M\subset \mathbb{R}^N$ such that $M$ is the closure of an open set and $M$ has smooth boundary. 
 Let $G=\mathbb{Z}$, which is amenable,  and consider the action $\beta: \mathbb{Z} \times M \rightarrow M$ given by $\beta_n = \varphi^n$, $n\in \mathbb{Z}$, with $\varphi$ a diffeomorphism of $M$ that fixes a non-empty open set of $M$. 
 As in Example \ref{example1}, consider the action of $\mathbb{Z}$ on $C(M)$ given by $\alpha_n(a)= a\circ \beta_n$, $a\in \cA$, $n\in \mathbb{Z}$. Clearly, $\alpha$  is not topologically free.

We will construct a covariant representation of $(\mathcal{A}, \mathbb{Z}, \alpha)$ that satisfies condition (B2). Consider the Hilbert space given by $H = L^2(M\times S^1)$, where $M$ has the usual Lebesgue measure and $S^1$ is endowed with the normalized Lebesgue measure, so that it has measure $1$. Consider the faithful representation $\pi:\mathcal{A} \rightarrow B(H)$ and the unitary representation $U : \mathbb{Z} \rightarrow B(H)$ given by
    \begin{equation*}
        [\pi(a)f](x,t)
    =
        a(x)
        f(x,t),
    \quad\quad
        [U_nf](x,t)
    =
        e^{int}
        |\det d\varphi^n|^{1/2}
        f(\varphi^{n}(x),t).
    \end{equation*}
    We thus obtain
    \begin{equation*}
        U_n\pi(a)U_n^*
    =
        \pi(\alpha_n(a)),
    \end{equation*}
    and so $(\pi, U)$ is a covariant representation of $(\cA, \mathbb{Z}, \alpha)$.

Let $G_0\subset \mathbb{Z}$ be finite and $\emptyset\neq V\subset M$ be any open set such that $\varphi^n|_V = id|_V$, for all $n\in G_0$. 
   Let $z_V: M \rightarrow \mathbb{R}$ be a continuous non-negative function  with support in $V$ and $f \in L^2(M\times S^1)$ be a function constant in the variable $t\in S^1$ such that $\|f\|_{L^2(M\times S^1)} = 1$. 
   Since $\varphi^n|_V = id|_V$, we have that $z_V(x)  |\det d\varphi^n|^{1/2} f(\varphi^n(x)) = z_V(x) f(x)$, thus for any $a_n\in \cA$, $n\in G_0$, 
    \begin{align*}
        &\left\|
         \sum_{n\in G_0}
            \pi(z_V)
            \pi(a_n)
            U_n
       \right \|_{B(L^2(M\times S^1))}^2 \;
         \geq \;
        \left \| \sum_{n\in G_0}
           \left [  \pi(z_V)
            \pi(a_n)
            U_n\right](f)
       \right \|_{L^2(M\times S^1)}^2\\
    = &
        \int_M
        \int_{S^1}
        \left|
            \sum_{n\in G_0}
            z_V(x)
            e^{int}
            a_n(x)
            |\det d\varphi^n|^{1/2}
            f(\varphi^{n}(x))
        \right|^2
        dt
        dx\\
    =&
        \int_M
        \int_{S^1}
        \left|
            \sum_{n\in G_0}
            z_V(x)
            a_n(x)
            e^{int}
            f(x)
        \right|^2
        dt
        dx\\
    =&
            \int_M
            \sum_{n\in G_0}   |z_V(x)|^2
           | a_n(x)|^2
            |f(x)|^2
            dx.
    \end{align*}
    Let $x_0 \in V$ be the maximum of the function $|z_V a_0|$.
    Now take a sequence $f_k:M \rightarrow \mathbb{R}$ such that $\|f_k\|_{L^2(M\times S^1)} = 1$ and $f_k^2$ is such that  for any  function $h\in C(M)$, we have $\int_M h(x) f_k^2(x) dx \to h(x_0)$, as $k\to \infty$ (for instance $f_k(x)= \chi_{ B_{1/k}(x_0)\cap M } (x)/ \sqrt{ m( B_{1/k}(x_0) \cap M)}$, with $m$ the Lebesgue measure and $B_{1/k}(x_0)$ the ball centered in $x_0$ with radius $1/k$). 
     Thus we obtain that
    \begin{align*}
         \lim_{k\to \infty}  \left  \|
       \sum_{n\in G_0}
          [  \pi(z_V)
            \pi(a_n)
            U_n](f_k)
       \right \|_{L^2(M\times S^1)}^2
    &=
        \sum_{n\in G_0} |z_V(x_0)|^2
       | a_n(x_0)|^2\\
    &\geq
     |   z_V(x_0)a_0(x_0)|^2.
    \end{align*}
We conclude that   
$$\left\|
         \sum_{n\in G_0}
            \pi(z_V)
            \pi(a_n)
            U_n
       \right \|_{B(L^2(M\times S^1))} \;
       \geq
     |   z_V(x_0)a_0(x_0)| =  
      \|z_Va_0\|,$$
  which proves condition (B2) in $\cB:= \alg(\pi(\cA), U_\mathbb{Z})$.  
   It then follows from Corollary \ref{cor.B2equivB1}  that the integrated representation $\pi\rtimes U$ gives an isomorphism of $\cA\rtimes_\alpha \mathbb{Z}$ with $\cB= \alg(\pi(\cA), U_\mathbb{Z})$. 

In this example, $\mathbb{Z}$ may be replaced by a generic discrete, amenable, group $G$, replacing $S^1$ by the compact space $\hat{G}$, which is the Pontryagin dual of $G$.
\end{example}


\subsection{Subalgebras and the commutative case}\label{ss.subalgebras_comm}

We consider here the equivalent conditions (B1) and (B2) on subalgebras of $\cB= \alg(\cA, U_G)$. We start with considering the class of operators with scalar coefficients, which, as we shall see, will be relevant in the commutative case.

 Let {$\tilde\cB := \alg(U_G)\subset \cB$} be the $C^*$-subalgebra generated by $U_G$, which is given by the closure of
 $$ \tilde\cB_0:=\tilde\cB \cap\cB_0 = \left\{\sum_{g \in G_{0}} (c_{g}I) U_{g}:  G_0 \; \text{finite},\; c_{g}\in \mathbb{C}\right \}.$$ 
 Let $b_{0}= \sum_{g \in G_{0}} (c_{g}I) U_{g}$ with $c_{g}\in \mathbb{C}$, then we can write
$$
b_{0}^{*}b_{0}  = \tilde{a_{0}}I +\sum_{s \neq t \in G_{0}} \overline{c_{s}} c_{t} U_{s^{-1}t} , \quad \text{with} \quad  \tilde{a_{0}}:= \sum_{s \in G_{0}} |c_{s}|^{2} \in \mathbb{C}
$$
Assuming $\tilde{a_{0}}>0$,  
 condition (B1) for $b_{0}$ can be written in the following way:  for any non-empty open set $V\subset M$, there exists $z_V\in \cZ(V)$ such that

$$\left\|z_{V}\left( \tilde{a_{0}} + \sum_{s \neq t \in G_{0}} \overline{c_{s}} c_{t}  U_{s^{-1}t}\right) \right\| 
= \tilde{a_{0}} \left\|z_{V}\left( I+ \sum_{s \neq t \in G_{0}} \overline{c_{s}} c_{t} /\tilde{a_{0}}  U_{s^{-1}t}\right) \right\| 
\geq \| z_{V} \tilde{a}_0\| = \tilde{a_{0}},  $$
that is, condition (B1) holds in $\tilde\cB\subset \cB$ if, and only if, for any finite $G_0\subset G$ and $c_{g} \in \mathbb{C}, g\in G_0$ and open $V \subset M$, there exists  $z_{V} \in \mathcal{Z}(V)$ such that
\begin{equation}\label{eq.B1'}
 \left\|z_{V}\left(I +\sum_{g \in G_0} c_{g} U_{g}\right)\right\| \geq\left\|z_{V}\right\|=1.
 \end{equation}

We can write condition (B2) in a similar way.
It follows from  Corollary \ref{prop.B2} that (B1) being satisfied in $\tilde\cB$ is equivalent to (B2)  being satisfied in $\tilde\cB$.

Assume  now $\cA$ is commutative. The point of considering the subalgebra $\tilde\cB$ is that in this case any $b\in \cB_0$ can be approximated by some {$b_{0}\in \tilde\cB_0=\tilde\cB \cap\cB_0$} on open sets.

\begin{proposition}\label{prop.locdense comm}
Assume $\cA=C(M)$  is commutative, $M$ a compact Hausdorff space. Let $b=\sum_{g\in G_0}a_gU_g$, with $a_g\in \cA$, and $\tilde a= \sum_{g\in G_0}\alpha_g^{-1}(a^*_ga_g)$. Given $\eps>0$, for any non-empty open $V\subset M$,
 and $m\in V$, there exists $b_{0}\in \tilde\cB_0$  and a non-empty open $\Delta \subset V$, $m\in \Delta$, such that 
$$\|z_{\Delta}( b^{*}b - b_{0}^{*}b_{0})\| < \eps \quad \text{and} \quad \|z_{\Delta}(  \tilde{a}-\tilde{a}_{0})\| <\eps $$
for any $z_{\Delta}\in  \cZ(\Delta)$ 
with $z_\Delta(m)=1$.
\end{proposition}

\begin{proof}
Let $\emptyset\neq V \subset M$ be open and fix $m \in V$. We have $b^*b=  \sum_{s,t\in G_0}\alpha_{s}^{-1}\left(a_{s}^{*} a_{t}\right)U_{s^{-1} t}$.  
  By continuity of the functions $a_s$,  we can take a neighborhood  $m \in \Delta \subset V$ small enough such that 
 $$\left|\alpha_{s}^{-1}\left(a_{s}^{*} a_{t}\right)(x)-\alpha_{s}^{-1}\left(a_{s}^{*} a_{t}\right)(m)\right|<\eps, \quad \mbox{ for all  }  x \in \Delta, \; s, t \in G_{0}$$
   Let $c_{s, t}=\alpha_{s}^{-1}\left(a_{s}^{*} a_{t}\right)(m) \in \mathbb{C}$, then for any $z_{\Delta}\in \cZ(\Delta)$, 
   assume also $z_{\Delta}(m)=1$, 
$$
\left\|z_{\Delta}\left(\alpha_{s}^{-1}\left(a_{s}^{*} a_{t}\right)-c_{s, t} I\right)\right\| 
= \sup_{x\in \Delta} z_\Delta(x)  \left|\alpha_{s}^{-1}\left(a_{s}^{*} a_{t}\right)(x) -c_{s, t} \right|
\leq \eps
$$
 Let $b_{0}= \sum_{g\in G_0}  a_g(m) U_{g}$ then $b_{0}^*b_0= \sum_{s,t\in G_0} c_{s,t} U_{s^{-1} t}$, and 
$$\|z_{\Delta}( b^{*}b - b_{0}^{*}b_{0})\| <K\eps, \qquad \|z_{\Delta}(  \tilde{a}-\tilde{a}_{0})\| <K'\eps
$$
for some $K, K'\in \N$ so the claim follows.
\end{proof}

Recall that we assume throughout  that we have a covariant representation $(id, U)$ of $(\cA, G, \alpha)$ with $\alpha$ an action that leaves a central subalgebra $\cZ$ invariant, that is, that condition (A1) holds. 

\begin{theorem}\label{prop.comm_red}  Let $\mathcal{A}=C(M)$ be commutative and assume (A1) is satisfied. If 
 (B2) holds in $\tilde \cB$,  that is, if
 for any finite $G_0\subset G$ and $c_{g} \in \mathbb{C}, g\in G_0$ and
  any non-empty open $V$ satisfying $\left.\beta_{g}\right|_{V}=\left.i d\right|_{V}$ for all $g \in G_0$, there exists  $z_{V} \in \mathcal{Z}(V)$ such that
\begin{equation}\label{eq.thm1}
 \left\|z_{V}\left(I +\sum_{g \in G_0} c_{g} U_{g}\right)\right\| \geq\left\|z_{V}\right\|=1,
 \end{equation}
then 
  (B2) holds in $\cB$. Assuming the action $\alpha$ is amenable, then  $id\rtimes U: \cA\rtimes_\alpha G \to \cB$ is an isomorphism. 
\end{theorem}

\begin{proof} 
 Let $b=\sum_{g\in G_0}a_gU_g\in \cB_0$,  $a_g\in \cA$, be arbitrary and $ \tilde a$ be as before, then from the previous proposition,  given $\eps>0$, for any open $V\subset M$ and $m\in V$, there exists $b_{0}\in \cB_0 \cap \tilde\cB$ and an open $\Delta \subset V$ such that 
$$\|z_{\Delta}( b_{0}^{*}b_{0})\| < \|z_{\Delta}( b^{*}b )\|  + \eps \quad \text{and} \quad \|z_{\Delta}(  \tilde{a}_{0})\| > \|z_{\Delta}(  \tilde{a})\| -\eps $$
for any $z_{\Delta}\in \cZ(\Delta)$.
  Assuming condition  \eqref{eq.thm1} gives that for any such open $V$ with  $\left.\beta_{g}\right|_{V}=\left.i d\right|_{V}$, 
 there is such a $z_{\Delta}$ satisfying
$\| z_{\Delta} b_{0}^{*}b_{0} \| \geq \| z_\Delta \tilde{a_{0}}\|$
hence $$\| z_{\Delta} b^{*}b \| \geq \| z_\Delta \tilde{a}\| - 2\eps$$ so that (B1) holds in $b$, hence (B2) holds in $\cB$. Since $\cA$ is commutative, it satisfies (C), hence from Corollary \ref{cor.B1thenB0}, 
(B0) holds and $id\rtimes U: \cA\rtimes_\alpha G \to \cB$ is an isomorphism.
\end{proof}

Hence, in the commutative case, we have that condition (B1), respectively, condition (B2), in $\mathcal{B}=\alg\left(\mathcal{A}, U_{G}\right)$ is equivalent to  condition (B1), respectively, condition (B2), in $\tilde\maB = \alg(U_G)$. Note that this condition depends on $M$-localization  and on the unitary group  $U_G$. 
We will see now that under an extra condition, the equivalent conditions (B1) and (B2) are  guaranteed just from conditions on $U_G$.
 
\begin{lemma}
    \label{pureStatesOverOpenLemma}
    Let $\tilde{A} \subset B(H)$ and $\tilde{B} \subset B(H)$ be commutative  unital  $C^*$-subalgebras that commute with each other, that is, for all $a \in \tilde A$ and $b \in \tilde B$ we have
    $ ab = ba.$
    Let $\tilde{A} \cong C(X)$ and $\tilde{B} \cong C(Y)$ and $\mathcal{C} = \text{alg}(\tilde{A}, \tilde{B}) \cong C(Z)$ and consider the induced maps by inclusion:
    $$
        \pi_X : Z \rightarrow X, \qquad
        \pi_Y : Z \rightarrow Y.
    $$
    If for all $0 \neq a \in \tilde{A}$ and $0 \neq b \in \tilde{B}$ we have
    \begin{equation}
    \label{nonTrivialMultiplication}
        ab \neq 0,
    \end{equation}
      then      for every neighborhood $x\in V \subset X$ we have that
    $
        \pi_Y(\pi_X^{-1}(V)) = Y.
    $
\end{lemma}
\begin{proof}
  Suppose  by contradiction that there exists $V \subset X$ open, $x\in V$, such that 
    $
        \pi_Y(\pi_X^{-1}(V)) \neq Y.
    $
    Since $X$ is normal (being compact and Hausdorff), we consider an open set $x \in U \subset V$ such that $\overline{U} \subset V$. By Uryshon's lemma, there exists $\rho \in \tilde{A} \cong C(X)$ such that:
    $$
        \rho(x) = 1, \quad \mbox{    and } \quad  \rho|_{X-\overline{U}} = 0.
    $$
    We have that $K=\pi_Y(\pi_X^{-1}(\overline U))\subset Y$ 
      is compact, thus $Y-K$ is a non-empty open set. We consider $y \in Y-K$, and by Uryshon's lemma, a function $g \in \tilde{B} \cong C(Y)$ such that
     $$
        g(y) = 1, \quad \mbox{
     and } \quad 
        g|_K = 0.
     $$
     Then we have that $ \rho g = 0 $:  if $z \in Z$ and $\pi_X(z) \notin \overline U$ then $\rho(z) = 0$,
    on the other hand, if $\pi_X(z) \in \overline U$, then $\pi_Y(z) \in K$ and so $g(z) = 0$.
    We obtain that $\rho(z)g(z) = 0$ for every $z \in Z$, thus $\rho g = 0$. This is in contradiction with \eqref{nonTrivialMultiplication}, since $\rho \neq 0$ and $g \neq 0$. 
 We have then that
    $
        \pi_Y(\pi_X^{-1}(V))
    =
        Y,
    $
     for all open sets $x\in V$ of $X$, which concludes the proof.
\end{proof}

We now have the following sufficient conditions for the isomorphism with the crossed product in the commutative case. Note that condition \eqref{normCondInGroupAlgebra} below is in fact (B0) for $\tilde\cB=\alg(U_G)$, so that it is equivalent to $\tilde\cB\cong \C\rtimes_\alpha G = C^*(G)$, the group algebra. It can be regarded as a strong form of linear independence of $U_g$, $g\in G$. 

\begin{theorem}\label{thm.condscomm2}
    Assume (A1) is satisfied and that $\mathcal{A} $ is commutative. Let $\tilde{\mathcal{B}} := \alg\{U_g: g\in G\} \subset B(H)$, 
    such that for all finite $G_0\subset G$, where $c_g \in \mathbb{C}$, we have
    \begin{equation}
    \label{normCondInGroupAlgebra}
        \left\|
            \sum_{g \in G_0}
            c_g
            U_g
        \right\|
        \geq
        |c_e|.
    \end{equation}
  Assume that for all $0 \neq a \in \mathcal{A}$ and $0 \neq b \in \tilde{\mathcal{B}}$ we have
    \begin{equation}
    \label{nonTrivialMultEq}
        ab \neq 0,
    \end{equation}
    then $(B2)$ holds in $\cB=\alg (\cA, U_G)$. If the action $\alpha$ is amenable, then $id\rtimes U: \cA\rtimes_\alpha G \to \cB$ is an isomorphism.
\end{theorem}

\begin{proof}

    By Theorem \ref{prop.comm_red}, it suffices to prove that given a finite set $G_0 \subset G$ and an open set $V \subset M$ such that the action $\beta_g|_V = id|_V$ is trivial for all $g \in G_0$, and a function $z_V \in \mathcal{Z}(V) \subset \mathcal{A}$ we have
  
    \begin{equation}
        \label{normConfFromThm}
        \left\|
            z_V
                \sum_{g \in G_0}
                c_g U_g
        \right\|
        \geq
        \|c_e z_V\|= |c_e|.
    \end{equation}
     Moreover,  we only need to prove \eqref{normConfFromThm} for positive  elements.
  Consider then $u = \sum_{g \in G_0} c_g U_g \in \tilde{\mathcal{B}}$ a positive element and $z_V \in \mathcal{Z}(V)$ such that $\beta_g|_V = id|_V$ for all $g \in G_0$. 
  
    From Remark \ref{rmk.opensets},  $u$ and $z_V$ commute. 
    We now consider $\overline{A} = \text{alg}\{z_V,Id\}$ and $\overline{B} = \text{alg}\{u, Id\}$. We have that $\overline{A}$ is a commutative $C^*$ algebra, since it is a sub-algebra of a commutative $C^*$ algebra and $\overline{B}$ is a commutative $C^*$ algebra since it is generated by a positive element $u$. Thus we have that $\overline{A} \cong C(X)$ and $\overline{B} \cong C(Y)$. Since $z_V$ commutes with $u$ we have that $\mathcal{C} := \text{alg}\{\overline{A}, \overline{B}\} \cong C(Z)$ is commutative, and by condition \eqref{nonTrivialMultEq}, since $\overline{A} \subset \mathcal{A}$ and $\overline{B} \subset \tilde{\mathcal{B}}$,  for all $0 \neq a \in \overline{A}$ and $b \in \overline{B}$ we have $ab \neq 0$. Thus we can apply Lemma \ref{pureStatesOverOpenLemma}. We consider  $x \in V$ such that
    $$
        z_V(x) = 1.
    $$
    We also consider $y \in Y$ such that
    $$
        u(y) = \|u\| \geq |c_e|.
    $$
    by equation \eqref{normCondInGroupAlgebra}. For $\eps>0$, consider the open set $V=\{x'\in X: z_V(x')>1-\eps\}$. Thus by Lemma \ref{pureStatesOverOpenLemma}, there exists $w \in Z$ such that $\pi_X(w) \in V $ and $\pi_Y(w) = y$. Then, with $x'=\pi_X(w)$,
    $$
        (z_V u)(w)
    =
        z_V(x')u(y)
    \geq
       (1-\eps) |c_e|.
    $$
    But then, since $w$ identifies with  a pure state in $\cC$, we conclude that
    $$
        \|z_Vu\|
    \geq
        (z_Vu)(w)
    \geq
        (1-\eps)|c_e|, \mbox{   for all  } \eps>0,
    $$
    which yields what we wanted to prove.
\end{proof}

 Put in a general framework where $\cA$ is an arbitrary unital $C^*$-algebra, not necessarily commutative, the results in this section can be regarded as density results in the following way.
 
Let  $\cA'$ be a $C^*$-subalgebra of $\cA$, containing $1_\cA$,  
  where $\alpha$ also acts by $\Aut(\cA')$ and let $\cB':=\alg(\cA', U_G)\subset \cB$, with $\cB_0'$ the set of elements of the form $\sum_{g\in G_0} c_gU_g$, $G_0$ finite, $c_g\in \cA'$. The conditions in Proposition \ref{prop.locdense comm} can be easily formulated in this setting:

For $b=\sum_{g\in G_0} a_gU_g\in \cB_0$ and $\tilde a=\sum_{g\in G_0} \alpha_g^{-1} (a_g^*a_g)\in \cA$, given $\eps>0$, and a non-empty open set $V\subset M$, 
$m\in V$, there exists $b_{0}=\sum_{g\in G_0} c_gU_g\in \cB'_0$  and a non-empty open $\Delta \subset V$, such that 
\begin{equation}\label{cond.locdense}
\|z_{\Delta}( b^{*}b - b_{0}^{*}b_{0})\| < \eps \quad \text{and} \quad \|z_{\Delta}(  \tilde{a}-\tilde{a}_{0})\| <\eps 
\end{equation}
for any $z_{\Delta}\in  \cZ(\Delta)$  with $z_\Delta(m)=1$, where $\tilde a_0=\sum_{g\in G_0} \alpha_g^{-1} (c_g^*c_g)\in \cA'$.

If  condition \eqref{cond.locdense} is  satisfied for any $b\in \cB_0$, we say that $\cB'$ and $\cA'$ are \emph{$M$-locally dense} in $\cB$ and $\cA$, respectively.
 We have shown that when $\cA$ is commutative,  then $\cA'=\C$ and $\cB'=\alg(U_G)$ are locally dense in $\cA$ and $\cB$, even though they are not dense with the strong operator topology.

Of course, in this situation, the proof of Theorem \ref{prop.comm_red} stands exactly in the same way, assuming now that (B1) holds in $\cB'$. 
(Note that if $\cA'$ contains $\cZ$, then this comes down to $\cB'\cong \cA'\rtimes_\alpha G$.) 

\begin{proposition}
Assume (A1) is satisfied. Let  $\cA'$ be a $C^*$-subalgebra of $\cA$ containing the identity,
  where $\alpha$ also acts by $\Aut(\cA')$, and let $\cB':=\alg(\cA', U_G)\subset \cB$, such that  $\cB'$ and $\cA'$ are $M$-locally dense subalgebras of $\cB$ and $\cA$, respectively.
If (B2) holds in $\cB'$,
then 
 (B2) holds in $\cB$. If (C) holds and the action $\alpha$ is amenable, then $id\rtimes U: \cA\rtimes_\alpha G \to \cB $ is an isomorphism. 
\end{proposition}

As for Theorem \ref{thm.condscomm2}, we see from the proof that, even in the non-commutative case,  if (B0) holds in $\cB'$, so in particular, assuming the action is amenable, $\cB'\cong \cA'\rtimes_\alpha G$, and if $ab\neq 0$ for $a\in \cA$, $b\in \cB'$, then (B1) also holds in $\cB'$. We then obtain:
\begin{theorem}\label{thm.condslocadense2}
Assume (A1) is satisfied and that the action $\alpha$ is amenable. Let  $\cA'$ be a $C^*$-subalgebra of $\cA$ containing the identity,
  where $\alpha$ also acts ameanably by $\Aut(\cA')$, and let $\cB':=\alg(\cA', U_G)\subset \cB$. Assume that $\cB'$ and $\cA'$ are $M$-locally dense subalgebras of $\cB$ and $\cA$, respectively, and that $(B0)$ holds in $\cB'$,
 in particular,  $\cB'\cong \cA'\rtimes_\alpha G$. 
 
 Assume that  $ab\neq 0$ for $0\neq a\in \cA$, $0\neq b\in \cB'$. 
Then $(B2)$ holds in $\cB$, and if (C) holds for $\cA$, then $id\rtimes U: \cA\rtimes_\alpha G \to \cB $ is an isomorphism.
\end{theorem}

\begin{proof}
Similarly to the proof of Theorem \ref{thm.condscomm2}, from Corollary \ref{cor.B2equivB1} and the previous proposition,  it suffices to check condition (B2) on $\cB'$. Given a finite set $G_0 \subset G$ and an open set $V \subset M$ such that the action $\beta_g|_V = id|_V$ is trivial for all $g \in G_0$, and a function $z_V \in \mathcal{Z}(V) \subset \mathcal{A}$, we want to show that, for any positive element $u = \sum_{g \in G_0} c_g U_g \in {\mathcal{B}_0'}$, $c_g\in \cA'$, 
  \begin{equation*}
        \left\|
            z_V
                \sum_{g \in G_0}
                c_g U_g
        \right\|
        \geq
        \|z_Vc_e \|.
    \end{equation*}
Again from  Remark \ref{rmk.opensets},  $u$ and $z_V$ commute, so that we can show, using Lemma \ref{pureStatesOverOpenLemma} and that, by assumption, $\|u\|\geq \|c_e\|$, that  
  $$
        \|z_Vu\|
    \geq
        (1-\eps)\|c_e\|, \mbox{   for all  } \eps>0,
    $$
    which shows that $     \|z_Vu\|\geq \|c_e\|\geq \|z_V c_e\|$. In particular (B2) holds in $\cB'$ and the result follows.
\end{proof}


\section{Back to the local trajectories method}\label{ss.localtraj_cond}

We consider here, as in Section \ref{s.localtraj}, the local trajectories representations and give alternative conditions to  establish the local trajectory method, Theorem \ref{thm.mainlocaltrajec}, in the  algebra $\B= \alg(\cA, U_{G})$. We assume, as always, that condition (A1) is verified, that is, that $(id, U)$ is a covariant representation of $(\cA, G, \alpha)$ preserving a central subalgebra $\cZ=C(M)$. We will replace condition (A2) that the group $G$ is amenable by the more general notion of amenability of the action $\alpha$. 

Let $\Omega$ be the orbit space of the induced action $\beta$ on $M$. Recall that for $\omega \in \Omega$,  $\pi_\omega$ is a representation on $B(\ell^2(G, H_\omega))$ such that for $a\in \cA$, $\xi\in \ell^2(G, H_\omega)$, $g, s\in G$, 
$$\left[\pi_\omega(a)\xi\right](s)
= \pi'_\omega( \alpha_{s}^{-1}(a))(\xi(s)), \qquad 
\left[\pi_{\omega}\left(U_{g}\right) \xi\right](s)=\lambda_g(\xi)(s)=\xi\left(g^{-1} s\right),  $$
 where $\pi'_\omega: \cA \to B(H_\omega)$ is $\pi'_\omega = \phi_\omega \circ \rho_\omega$ with  $\rho_\omega= \rho_{m_\omega}: \cA\to \cA/J_{m_\omega}$ the quotient map, 
and  $\phi_\omega$ an isometry. 
 
 Then $\pi_\omega$ is defined in $\cA\rtimes_\alpha G$ as the regular representation induced by $\pi_\omega'$. As for $\cB=\alg(\cA, U_G)$, we assume that elements $b= \sum_{g\in G_0}a_gU_g$ uniquely determine the coefficients $a_g$, $g\in G_0$, so that $\pi_\omega$ is also well-defined in  $\cB_0$. 
Typically, we are interested in the case when $\cB\cong \cA\rtimes_\alpha G$ so this condition is guaranteed.

Throughout this section, we let $\pi = \oplus_{\omega \in \Omega} \pi_\omega$. As a matter of terminology, we say that  the local trajectories method works on $\cB$ if $\pi$ is well-defined and faithful in $\cB$. Recall from Proposition \ref{prop.injectiveA} that $\pi$ is always faithful in $\cA$. 

We have noted in Remark \ref{rmk_isom_loctraj} that for the local trajectories method to work on $\cB$ it suffices to give conditions such that  the maps
\begin{equation}\label{eq.isomlocaltraj}
i d \rtimes_\alpha U: \mathcal{A} \rtimes_\alpha G \rightarrow \mathcal{B}, \qquad
 i d \rtimes_{\alpha^{\prime}} \pi(U): \pi(\mathcal{A}) \rtimes_{\alpha^{\prime}} G \rightarrow \pi(\mathcal{B})
\end{equation}
are isomorphisms, where   $\alpha'_g(\pi(a)):= \pi(U_g)\pi(a) \pi(U_g^*) = \pi( U_ga U_g^*)$. Note that $\pi(\cB)\cong \alg(\pi(\cA), \pi(U_G))$.

The goal of this section is to apply the criteria obtained in the Section \ref{s.cond_isom} to study in particular the second map in \eqref{eq.isomlocaltraj}, where we have regular representations on spaces of the form $B(\ell^2(G, H))$. 

We consider first  the commutative case where we can use the results in Section \ref{ss.subalgebras_comm}. Let $ \tilde\cB_0 = \tilde\cB\cap \cB_0$ and $\tilde \cB= \alg(U_G)$ and consider operators of the form
 $$b_0= I +\sum_{e\neq s\in G_0} c_sU_s\; \in\;  \tilde\cB_0, \quad  \mbox{with    } \pi_\omega(b_0)=I+\sum_{e\neq s\in G_0} c_s\lambda_s, c_s\in \C.$$ 
 Then for $\xi\in \ell^2(G, H_\omega)$ and $g\in G$, we have
 $$\left[\pi_\omega(b_0)(\xi)\right](g)= \xi(g) + \sum_{g\neq s\in G_0} c_s\xi(g^{-1}s).$$
   We have that  the localized condition (B1) always holds in $\pi(\tilde\cB_0)$. 
   For that we make use of the following lemma.

\begin{lemma}\label{lem.rho=1}
Let  $\rho\in \cZ\cong C(M)$ such that $0\leq \rho\leq 1$ and $\| \rho\|=1$. Then there exist $\omega \in \Omega$ and $g\in G$ such that 
\begin{equation}\label{eqident}
\alpha_g^{-1}\left(\rho\right) = 1_\cA + J_{m_\omega}.
\end{equation}
In particular, $\pi'_\omega(\alpha_g^{-1}\left(\rho\right) ) =I\in B(H_\omega)$ and 
$[\pi_\omega(\rho)\xi] (g)= \xi(g)$, for any $\xi\in  \ell^2(G, H_\omega).$
Conversely, given $\omega \in \Omega$ and $g\in G$, there is $\rho\in \cZ\cong C(M)$ such that $0\leq \rho\leq 1$, and $\| \rho\|=1$ satisfying \eqref{eqident}.
\end{lemma}

\begin{proof}
Let $m\in M$ be such that $\rho(m)=1 = \max_{x\in M} \rho(x)$.
Take the orbit $\omega \in \Omega$ such that $m \in \omega$,   so  $m=\beta_g^{-1}\left(m_\omega\right)$ for $m_\omega \in \omega$ and $g \in G$. Note that 
$\alpha_g^{-1}\left(\rho\right)=\rho \circ \beta_g^{-1}$ will have its maximum  in $m_\omega$ since $\beta_g^{-1}\left(m_\omega\right)=m$ and $\rho$ attains its maximum at  $m$. This gives that $\alpha_g^{-1}(\rho)(m_\omega)=1$, hence $\alpha_g^{-1}(\rho)- 1_{\cA} \in J_{m_\omega}$ and
$$\pi'_\omega \left(\alpha_g^{-1}\left(\rho\right)\right)= \pi'_{m_\omega}\left(\alpha_g^{-1}\left(\rho\right)\right)= \pi'_{\omega}(1)= I\in B(H_\omega),$$ from which follows $\pi_\omega(\rho)\eta (g)= \pi'_\omega \left(\alpha_g^{-1}\left(\rho\right)\right)\eta(g)= \eta(g)$, for any $\eta\in  \ell^2(G, H_\omega).$
For the converse, let $\omega \in \Omega$ and $g\in G$ and pick any $\rho\in \cZ$ such that $0\leq \rho\leq 1$ and $ \rho(m)=1$, with $m= \beta_g^{-1}(m_\omega)$.
\end{proof}

\begin{proposition}\label{prop.B1scalar}
Assume condition (A1).  Let $\pi= \oplus_{\omega \in \Omega} \pi_\omega$ be the local trajectories representation, assumed well defined in $\cB_0$. Then  for any $V\subset M$ open and $c_s\in \C$, $s\in G_0\subset G$ finite, and for all $z_V\in \cZ(V)$ we have 
$$ \left\|z_V\pi\left(I +\sum_{e \neq s \in G} c_s U_s\right)\right\|\geq 1.$$
In particular, (B1) holds in  $\pi(\tilde\cB_0 )\subset \pi(\cB)$. 
\end{proposition}

\begin{proof}
Let $z_V\in \cZ(V)$, then since $\pi$ is an isomorphism on $\cZ$, write $z_V= \pi(\rho_V)$, for $\rho_V\in \cZ$. Then $\|\rho_V\|=\|\pi(\rho_V)\|=1$ and since $\pi$ is a $*$-homomorphism and $z_V\geq 0$, also $\rho_V\geq 0$. 

 Let $g\in G$ and $\omega\in \Omega$ be as in  the previous lemma applied to $\rho_V$.  Then we have 
  $$\left(\pi_\omega(\rho_V)\xi \right)(g)
  = \pi'_\omega \left(\alpha_g^{-1}\left(\rho_V\right)\right)\left(\xi(g)\right)
  = \xi(g),$$ for any $\xi\in  \ell^2(G, H_\omega)$.  
 Now take  $\xi_g \in \ell^2(G, H_\omega)$ such that  $\xi_g(s)=0$ for $s \neq g$ and $\xi_g(g)=u$, with $u\in H_\omega$, $\|u\|=1$, arbitrary. Let $b_0=I+ \sum_{e\neq s\in G_0} c_sU_s\in \tilde{\cB}_0$. It follows that
\begin{equation*}\label{eq.eqbetag}
\left(\pi_\omega\left(\rho_V b_0\right) \xi_g\right)(g)=
\pi_\omega(\rho_V)\left(\xi_g(g)+\sum_{e \neq s \in G} c_s \xi_g\left(s^{-1} g\right)\right)=\xi_g(g)=u.
\end{equation*} 
 We have then
 \begin{align*}\label{eq.normgreater1}
 \|\pi_\omega(\rho_V b_0)\| & = \sup_{\|\xi\|=1} \|\pi_\omega(\rho_V b_0)\xi\|_{\ell^2} 
  \geq \|\pi_\omega(\rho_V b_0)\xi_g\|_{\ell^2} \\
  &\geq  \|(\pi_\omega(\rho_V b_0)\xi_g)(g)\|=\|u\|=1.
 \end{align*}
 Since we wrote $z_V=\pi(\rho_V)$, we conclude that
$$ \left\|z_V\pi\left(I +\sum_{e \neq s \in G} c_s U_s\right)\right\|= \sup_{\omega \in \Omega} \left\|\pi_\omega(\rho_V)\pi_\omega\left(I +\sum_{e \neq s \in G} c_s U_s\right)\right\|\geq 1.$$
In particular, (B1) holds in $\pi(\tilde \cB_0)$ .
\end{proof}

Assume now that $\cA$ is commutative; then both $\cA$ and $\pi(\cA)$ satisfy the $M$-localization condition (C) as in \eqref{C}. 
Moreover, since $\tilde\cB_0=\cB_0 \cap \tilde\cB$ is $M$-locally dense in $\cB$, by Proposition \ref{prop.locdense comm}, also $\pi(\tilde\cB_0) = \pi(\cB_0) \cap \alg(\pi(U_G))$ is $M$-locally dense in $\pi(\cB)$,  so we have from Theorem \ref{prop.comm_red} that if (B2), or equivalently, (B1), holds in $\pi(\tilde\cB_0)$ then it also holds in $\pi(\cB)$. In this case, using localization, it follows from Corollary \ref{cor.B2equivB1} that:

\begin{proposition}
\label{lem_commB2}
 Assume condition (A1) and that the action $\alpha$ is amenable. Let $\cA$ be commutative with the local trajectories representation $\pi=\oplus_{\omega \in \Omega}\pi_\omega$ well defined in $\cB_0$. Then (B1) holds in $\pi(\cB)$ and
$$id\rtimes \pi(U): \pi(\cA)\rtimes_{\alpha'} G \to \pi(\cB)
 = \alg(\pi(\cA), \lambda_G)$$
 is an isomorphism.
\end{proposition}

Following Remark \ref{rmk_isom_loctraj}, we conclude that, in the commutative case, for the local trajectories method to work, it suffices  that $i d \rtimes_\alpha U: \cA \rtimes_\alpha G\to \cB$ is an isomorphism. From Theorem \ref{prop.comm_red} we obtain then the following sufficient condition.
\begin{theorem}\label{thm.localtraj_comm}
 Let $\cA$ be commutative. 
 Assume condition (A1) and that the action $\alpha$ is amenable. 
  If 
  (B2) holds in the subalgebra $\tilde\cB_0$, that is, if
 for any finite $G_0\subset G$ and $c_{g} \in \mathbb{C}, g\in G_0$ and 
  any non-empty open set $V$ satisfying $\left.\beta_{g}\right|_{V}=\left.i d\right|_{V}$ for all $g \in G_0$, there exists  $z_{V} \in \mathcal{Z}(V)$ such that
\begin{equation*}\label{eq.thm}
 \left\|z_{V}\left(I +\sum_{g \in G_0} c_{g} U_{g}\right)\right\| \geq\left\|z_{V}\right\|=1,
 \end{equation*}
then $\{\pi_{\omega}\}_{\omega\in \Omega}$ is a faithful family of representations of $\cB$, that is, the local trajectories method works on $\cB$. 
\end{theorem}

 We now show that in fact  Proposition \ref{lem_commB2}  still holds 
 even if $\cA$ is not commutative, relying essentially on properties of regular representations. 

For elements $b= \sum_{s\in G_0} a_sU_s\in \cB_0$, we have 
\begin{align*}
\left(\pi_\omega\left( b\right) \xi\right) (g) & =  \sum_{s\in G_0} 
  \left( \pi'_\omega(\alpha_g^{-1}(  a_s) \right) (\xi( s^{-1}g ))
\end{align*}
and with  $\xi \in \ell^2(G, H_\omega)$ such that  $\xi(t)=0$ for $t \neq g$ and $\xi(g)=u$,  $u\in H$, $\|u\|=1$, 
\begin{align*}\label{eq.piomega}
\left(\pi_\omega\left( b\right) \xi\right) (g) & =  
 ( \pi_\omega( a_e)\xi)(g) = \left( \pi'_\omega(\alpha_g^{-1}(  a_e)\right) (\xi(g))
\end{align*}
Introducing  operators $j_g: H\to \ell^2(G, H)$ and $j_g^*: \ell^2(G, H)\to H$, $g\in G$, such that, for $h\in H$, $\xi \in  \ell^2(G, H)$, 
 \begin{equation}\label{eq.Jg}
 j_g(h)(g)= h,\; \; j_g(h)(t)= 0, t\neq g, \quad  \mbox{ and } \quad j_g^*(\xi) = \xi(g), 
  \end{equation}
  we have that $j_g$ is an isometry and $\|j_g\|=\|j_g^*\|=1$, with $j_g^* j_g=I$. 
   We can then write the equality above as 
 \begin{equation}\label{eq.Jg2}
 j_g^* \pi_\omega( b) j_g= \pi'_\omega(\alpha_g^{-1}( a_e)) \; \mbox{ in }B(H_\omega).
 \end{equation}

 We
make use of the following lemma (that holds in general for regular representations, similarly to the result just after).

\begin{lemma}\label{lem_normpia}
For each $\omega\in \Omega$, $a\in \cA$, 
$$\displaystyle \|\pi_\omega(a)\|
= \sup_{g\in G} \|\pi_\omega'(\alpha_g^{-1}(a))\|_{B(H_\omega)}.$$
\end{lemma}
\begin{proof}
For each $g\in G$, we can write,  as in \eqref{eq.Jg2}, 
$ j_g^* \pi_\omega(a) j_g= \pi'_\omega(\alpha_g^{-1}(a))$. Hence, since $\|j_g\|=\|j_g^*\|=1$, 
$$ \| \pi'_\omega(\alpha_g^{-1}(a))\|= \|j_g^* \pi_\omega(a) j_g \|\leq \| \pi_\omega(a) \|$$
It follows that  $\sup_{g\in G} \|\pi_\omega'(\alpha_g^{-1}(a))\|_{B(H_\omega)}\leq \| \pi_\omega(a) \|$.
For the reverse, let $\xi \in \ell^2(G, H_\omega)$ with $\|\xi\|^2=    
        \sum_{g \in G}
        \|\xi(g)\|^2_{H_w}\|
=1$.
    Then we have that:
    \begin{align*}
        \|\pi_w(a) \xi\|^2_{\ell^2(G,H_w)}
    &=
        \sum_{g \in G}
        \|\pi'_{\omega}(\alpha_g^{-1}(a))\xi(g)\|^2_{H_\omega} \leq
        \sum_{g \in G}
        \|\pi'_{\omega}(\alpha_g^{-1}(a))\|_{B(H_\omega)}\|\xi(g)\|^2_{H_w}\\
    &\leq
        \sup_{g \in G}\|\pi'_\omega(\alpha_g^{-1}(a))\|_{B(H_\omega)}
        \sum_{g \in G}
        \|\xi(g)\|^2_{H_\omega} \leq
        \sup_{g \in G}\|\pi'_\omega(\alpha_g^{-1}(a))\|_{B(H_\omega)}.
    \end{align*}
    Taking the supremum over $\xi \in \ell^2(G, H_\omega)$ such that $\|\xi\|_{\ell^2(G,H_\omega)} = 1$ we obtain:
    $$
        \|\pi_w(a)\|_{\ell^2(G,H_\omega)}
    \leq
        \sup_{g \in G}
        \|\pi_\omega'(\alpha_g^{-1}(a))\|_{B(H_\omega)},
    $$
    which concludes the proof.
\end{proof}

A first consequence is the following version of (B0) in $\pi_\omega(\cB)$.
\begin{proposition}\label{lemma.B0piomega}
Assume condition (A1) and that the local trajectories representation $\pi=\oplus_{\omega\in \Omega} \pi_\omega$ is well defined in $\cB_0$. Then, for any finite set $G_0\subset G$ and $b=\sum_{s\in G_0} a_sU_s$, and for any orbit $\omega\in \Omega$, we have
\begin{equation}\label{eq.B0piomega}
\|\pi_\omega(b)\|\geq \|\pi_\omega(a_e)\|.
\end{equation}
\end{proposition}
\begin{proof}
With the notation as in \eqref{eq.Jg}, we have for each $\omega \in \Omega$ and all $g\in G$, $ j_g^* \pi_\omega(b) j_g= \pi'_\omega(\alpha_g^{-1}( a_e))$, hence
$$\| \pi'_\omega(\alpha_g^{-1}( a_e))\|= \|  j_g^* \pi_\omega(b) j_g\|\leq \|\pi_\omega(b)\|, \quad \mbox{ for all } g\in G.$$
Since, from Lemma \ref{lem_normpia}, $\|\pi_\omega(a_e)\|= \sup_{g\in G} \| \pi'_\omega(\alpha_g^{-1}( a_e))\|$, the result follows.
\end{proof}

We  now show that the $M$-local condition (B1)  holds in $\pi(\cB)$ so the second isomorphism in \eqref{eq.isomlocaltraj} always  holds, in case our algebra $\cA$ is $M$-localizable, that is, satisfies (C).

\begin{proposition}\label{prop.extracond}
Assume condition (A1) and that the action $\alpha$ is amenable. Let $\pi= \oplus_{\omega \in \Omega} \pi_\omega$ be the local trajectories representation, assumed well defined in $\cB_0$. Then  for any $V\subset M$ open and $a_s\in \cA$, $s\in G_0\subset G$ finite, and for all $z_V\in \cZ(V)$ we have 
$$ \left\|z_V\pi\left(a_e+\sum_{e \neq s \in G_0} a_s U_s\right)\right\|\geq  \left\|z_V\pi\left(a_e\right)\right\|.$$
In particular, 
  (B1) holds in $\pi(\cB)$,  and if (C) holds for $\cA$, then
  $$id\rtimes \pi(U): \pi(\cA)\rtimes_{\alpha'} G \to \pi(\cB)
 $$
 is an isomorphism.
\end{proposition}
\begin{proof}
Similarly to the proof of Proposition \ref{prop.B1scalar}, writing $z_V= \pi(\rho_V)$, it suffices to prove that for any open $V\subset M$,  $\rho_V\in \cZ(V)$ and $b= \sum_{s\in G_0} a_sU_s\in \cB_0$, we have  $\|\pi(\rho_V b)\|\geq \|\pi(\rho_V a_e)\|$. This is straightforward from the previous proposition, since we have
$\|\pi_\omega(\rho_Vb)\|\geq \|\pi_\omega(\rho_Va_e)\|$, for all $\omega \in \Omega$.
 In particular,  (B1) holds in $\pi(\cB)$, and the isomorphism follows from $\pi(\cA$) also satisfying (C) and Corollary \ref{cor.B1thenB0}.
\end{proof}

It then follows again from Remark \ref{rmk_isom_loctraj} that, if $\cA$ is $M$-localizable, then the local trajectories method works on $\cB$, as long as  $i d \rtimes_\alpha U: \cA \rtimes_\alpha G\to \cB$ is an isomorphism. We then obtain the following $M$-localized version of the local trajectories method, Theorem \ref{thm.mainlocaltrajec}. (Note that if the action is topologically free, then there are no open sets of fixed points, so the condition is trivially verified.)

\begin{theorem}\label{thm.extracond}
Assume condition (A1) and that the action $\alpha$ is amenable. If (C) holds for $\cA$ and the $M$-local condition (B2) holds in $\cB$, that is, if for every  finite set $G_0 \subset G$ and non-empty open set $V \subset M$ such that $\left.\beta_{g}\right|_{V}=\left.i d\right|_{V}$ for all $g \in D$, and for  $a_{g} \in \cA$, there exists  $z_{V} \in \mathcal{Z}(V)$  
such that
 $$\left\|z_{V}\sum_{g  \in G_0} a_{g} U_{g}\right\| \geq\left\|z_{V}{a_e}\right\|,$$
 then $\{\pi_{\omega}\}_{\omega\in \Omega}$ is a faithful family of representations of $\cB$, that is, the local trajectories method works on $\cB$. 
\end{theorem}

Along the lines of the last results of the previous section, Theorems \ref{thm.condscomm2} and \ref{thm.condslocadense2},
we also conclude:

\begin{theorem}\label{thm.condscomm2.localtraj}
   Let $\mathcal{A} $ be commutative. Assume condition (A1) and that the action $\alpha$ is amenable. Let $\tilde{\mathcal{B}} := \alg\{U_g: g\in G\} \subset B(H)$, 
    and assume that for all finite $G_0\subset G$, where $c_g \in \mathbb{C}$, we have
    \begin{equation*}
    \label{normCondInGroupAlgebra2}
        \left\|
            \sum_{g \in G_0}
            c_g
            U_g
        \right\|
        \geq
        |c_e|,
    \end{equation*}
that is, $\tilde\cB\cong \C\rtimes_\alpha G = C^*(G)$, the group algebra.  Assume that for all $0 \neq a \in \mathcal{A}$ and $0 \neq b \in \tilde{\mathcal{B}}$ we have
    \begin{equation*}
    \label{nonTrivialMultEq}
        ab \neq 0.
    \end{equation*}
    Then 
     $ \cA\rtimes_\alpha G \cong \cB$ and $\{\pi_{\omega}\}_{\omega\in \Omega}$ is a faithful family of representations of $\cB$, that is,  the local trajectories method works on $\cB$.
\end{theorem}

In general, if we have an isomorphism with the crossed product on $M$-locally dense subalgebras, as in \eqref{cond.locdense}, then:

\begin{theorem}\label{thm.condslocadense2.localtrajec}
Assume (A1) is satisfied and that the action $\alpha$ is amenable. Let  $\cA'$ be a $C^*$-subalgebra of $\cA$ containing the identity,
  where $\alpha$ also acts ameanably by $\Aut(\cA')$, and let $\cB':=\alg(\cA', U_G)\subset \cB$. Assume that $\cB'$ and $\cA'$ are $M$-locally dense subalgebras of $\cB$ and $\cA$, respectively, and that $(B0)$ holds in $\cB'$,
 in particular,  $\cB'\cong \cA'\rtimes_\alpha G$.

  Assume that  $ab\neq 0$ for $0\neq a\in \cA$, $0\neq b\in \cB'$. 
Then 
 if (C) holds for $\cA$, we also have
 $\cB\cong \cA\rtimes_\alpha G$ 
and $\{\pi_{\omega}\}_{\omega\in \Omega}$ is a faithful family of representations of $\cB$, that is,  the local trajectories method works on $\cB$.
\end{theorem}

\bibliographystyle{plain}
\bibliography{/Users/catarinacarvalho/Documents/Suf_families_local_traj/localtrajbib}

The three authors were supported by FCT--Funda\c c\~ao para a Ci\^encia e Tecnologia, Portugal, 
under project UIDB/04721/2020. C.C was also supported by the FCT-Portugal project UIDB/04459/2020, with DOI identifier 10-54499/UIDP/04459/2020. M.D was also supported by Research Foundation– Flanders (FWO) via the Odysseus II programme no. G0DBZ23N.

\end{document}